\documentclass{siamart0216}

\usepackage{lipsum}
\usepackage{amsfonts}
\usepackage{graphicx}
\usepackage{epstopdf}
\usepackage{algorithm}
\usepackage[noend]{algpseudocode}
\usepackage{mwe,tikz}
\usetikzlibrary{positioning,calc,fadings,decorations.pathreplacing}
\usepackage{tikz-3dplot}
\usetikzlibrary{arrows, decorations.markings}
\usetikzlibrary{calc,shapes, positioning}
\usetikzlibrary{matrix}

\usepackage{tikz}
\usetikzlibrary{shapes,shadows,trees}
\ifpdf
  \DeclareGraphicsExtensions{.eps,.pdf,.png,.jpg}
\else
  \DeclareGraphicsExtensions{.eps}
\fi

\newcommand{\TheTitle}{Manycore parallel computing for a hybridizable discontinuous Galerkin nested multigrid method} 
\newcommand{\TheAuthors}{M. S. Fabien, M. G. Knepley, R. T. Mills, and B. M. Rivier\`{e}}

\headers{Manycore parallel computing for HDG GMG}{\TheAuthors}

\title{{\TheTitle}}

\author{
  Maurice S. Fabien\thanks{Department of Computational and Applied Mathematics, Rice University, Houston, TX.}
  \and
  Matthew G. Knepley\thanks{Department of Computer Science and Engineering, University of Buffalo, Buffalo, NY.}
    \and
  Richard T. Mills\thanks{Laboratory for Advanced Numerical Simulations, Argonne National Laboratory.}
 \and
  B\'{e}atrice M. Rivi\`{e}re$^*$
}

\usepackage{amsopn}

\usepackage{bm}
\usepackage{booktabs}
\usepackage{xcolor,colortbl}
\usepackage{enumerate} 
\usepackage{subfig}
\usepackage{amsmath}

\usepackage{tikz}

\usepackage{tikz-cd}
\usetikzlibrary{shapes,arrows,intersections}
\usetikzlibrary{matrix,fit,calc,trees,positioning,arrows,chains,shapes.geometric,shapes}
 
\usetikzlibrary{intersections}

\usepackage{pgfplots}
\pgfplotsset{compat=newest,compat/show suggested version=false}
\usepackage{pgfplotstable}
\usepackage{filecontents}

\definecolor{Gray}{gray}{0.85}
\definecolor{LightCyan}{rgb}{0.88,1,1}

\newcolumntype{a}{>{\columncolor{Gray}}c}
\newcolumntype{b}{>{\columncolor{white}}c}

\ifpdf
\hypersetup{
  pdftitle={\TheTitle},
  pdfauthor={\TheAuthors}
}
\fi

\begin{document}
\fontsize{10}{10}\selectfont
\maketitle

\begin{abstract}
  We present a parallel computing strategy for a hybridizable discontinuous Galerkin (HDG) nested geometric multigrid (GMG) solver.  Parallel GMG solvers require a combination of coarse-grain and fine-grain parallelism to improve time to solution performance.  In this work we focus on fine-grain parallelism.  We use Intel's second generation Xeon Phi (Knights Landing) many-core processor.  The GMG method achieves ideal convergence rates of $0.2$ or less, for high polynomial orders.  A matrix free (assembly free) technique is exploited to save considerable memory usage and increase arithmetic intensity.  HDG enables static condensation, and due to the discontinuous nature of the discretization, we developed a matrix vector multiply routine that does not require any costly synchronizations or barriers.  Our algorithm is able to attain 80\% of peak bandwidth performance for higher order polynomials.  This is possible due to the data locality inherent in the HDG method. Very high performance is realized for high order schemes, due to good arithmetic intensity, which declines as the order is reduced.     
\end{abstract}



\section{Introduction}
Multigrid methods are among the most efficient solvers for linear systems that arise from the discretization of partial differential equations.  The effectiveness of this multilevel technique was first observed by~Fedorenko (\cite{Fedorenko}) in 1964, and popularized by Brandt (\cite{Brandt1}) in 1977.  Traditionally, multigrid methods have been applied to low order finite difference, finite volume, and continuous finite element approximations discretizations (\cite{trottenberg}, \cite{wesseling}, \cite{hackbusch}, \cite{Brandt2}).  The discontinuous Galerkin (DG) discretization brings with it many advantages: can handle complex geometries, has access to $hp$-adaptivity, is capable of satisfying local mass balance (ideal for flow problems and hyperbolic PDEs), and highly parallelizable due to the lack of continuity constraints between elements.  However, DG methods often have more degrees of freedom than their continuous counterpart.  In addition, high order discretizations rapidly increase the condition number of the discretization operator, which poses a challenge for any linear solver (\cite{Bastian}).
\\
\indent  Originally multigrid methods were developed at a time when access to parallelism was limited, and iterative methods that had less concurrency but better convergence rates were favored.  The multilevel nature of multigrid can cause challenges in its parallelization; load balancing problems occur due to fine grids have ample data to work with, but coarser levels do not.  Moreover, traditional multigrid is a multiplicative method, that is, each level must be completely processed before moving to the next.  Additive multigrid methods allow for the simultaneous processing of all levels, but the trade off between concurrency and robustness is often not ideal \cite{bastian1998additive}.
\\
\indent  Multiplicative multigrid is well known to have sequential complexity, for $N$ data points a single cycle costs $\mathcal{O}(N)$ floating point operations.  However, the parallel complexities of V- and F-cycles are polylogarithmic.  For multiplicative multigrid, a natural heterogeneous computing strategy is to process fine levels up to a threshold on an coprocessor (or accelerator), and have the remaining coarse levels be processed by coarse-grain parallelism.  In this paper all computations are done on a single Xeon Phi coprocessor (Knights Landing), as we only focus on the fine-grain parallelism.  However, an offloading model is a natural extension of this work.
\\
\indent
In this paper we consider a nested multigrid technique for high order hybridizable discontinuous Galerkin methods; which combines $p$-multigrid and $h$-multigrid.  The HDG method is capable of static condensation, which significantly reduces the total number of degrees of freedom compared to classical DG methods.  This has important consequences for linear solvers like multigrid, since the cost of multigrid grows with the number of degrees of freedom.  Moreover, due to the discontinuous nature of HDG methods, we show that we can efficiently leverage the massive parallelism that many-core processors offer.  Through roofline performance modeling and tuning, we show that our implementation results in efficient device utilization as well as significant speedups over a serial implementation.  An interesting benefit of the Xeon Phi coprocessor is that it allows us to use traditional parallel programming paradigms like OpenMP, MPI, and pthreads.  Consequently, we can take advantage of the massive fine-grain parallelism that the Xeon Phi offers by utilizing these traditional parallel paradigms with significantly limited software intrusion.

The remainder of the paper is organized as follows.  Section~\ref{sec:model} introduces the model PDE.  The relevant finite element notations and HDG discretization is described in section~\ref{sec:disc0}.  In section~\ref{sec:basis}, we explain how we implement the nodal tensor product basis, quadrature, and barycentric interpolation.  Numerical experiments verifying the correct convergences rates for the HDG discretization is conducted in~section~\ref{sec:disc}.  The core components of the multigrid solver we use is discussed in section~\ref{sec:mg_description}.  This includes the description of intergrid transfer operators, relaxation, coarse grid operators, and the standard multigrid cycle.  In section~\ref{sec:gmg} we show through computational experiments that our multigrid method obtains very good convergence rates and error reduction properties.  The performance of our algorithm in a parallel setting is evaluated in section~\ref{sec:matvec}.  
 

\section{Model problem}
\label{sec:model}
Consider the elliptic problem
\begin{align}
-\nabla \cdot ({\bf K} \nabla   u) &= f, \quad \mbox{in}\quad \Omega, 
\label{eq:pde1}
\\
u &= g_D,\quad\mbox{on}\quad\partial\Omega.
\label{eq:pde2}
\end{align}
where $\Omega$ is an open domain in $\mathbb{R}^2$ and $\partial\Omega$ denotes the boundary of the domain.  Dirichlet datum $g_D$ is imposed on the boundary.  The function $f$ is the prescribed source function and the matrix $\bf K$ is symmetric positive definite with piecewise constant entries. The div-grad operator in \eqref{eq:pde1} appears in several problems in engineering such as multiphase flow in porous media.  It will provide insight into how a DG multigrid solver performs on modern architectures.  Moreover, it acts as a gateway to construct very efficient numerical methods for time dependent PDEs that require pressure solves or implicit time stepping.  

\section{Discretization}
\label{sec:disc0}
HDG methods were designed to address the concern that DG schemes generate more degrees of freedom (DOFs) when compared to continuous Galerkin techniques.  For standard DG methods, each degree of freedom is coupled with all other degrees of freedom on a neighboring element.  By introducing additional unknowns along element interfaces, the HDG method is able to eliminate all degrees of freedom that do not reside on the interfaces.  As such, a significantly smaller linear system is generated, and HDG gains much of its efficiency at higher orders (\cite{Kirby2012}).  It turns out that the HDG method also has a number of attractive properties, namely, the capability of efficient implementations, optimal convergence rates in the potential and flux variables, as well as the availability of a postprocessing technique that results in the superconvergence of the potential variable.  HDG methods are a subset of DG methods, so they still retain favorable properties, e.g., local mass balance, ease of $hp$-adaptivity, and the discontinuous nature of the solution variables.  A thorough analysis of HDG methods can be found in \cite{CockburnDG08}, \cite{CockburnGL09}, and \cite{CockburnDGRS09}.
\\
\indent  A number of works are available on multigrid for DG methods.  Interior penalty methods are the most commonly analyzed, for instance, see \cite{brenner2005convergence}, \cite{brenner2009multigrid}, \cite{gopalakrishnan2003multilevel}, \cite{Antonietti2016}, and \cite{Ant}.  Most of these works are theoretical, and while they are able to prove convergence, the numerical experiments show rates below what is typically expected from GMG in this model setting.  For the interior penalty class of DG methods, it was found that specialized smoothers and careful tuning of the stability parameter was required for better convergence results (see \cite{johannsen2005symmetric}, \cite{johannsen2005multigrid}).  Local Fourier analysis (local mode analysis, \cite{Brandt2}) was applied to interior penalty DG methods in \cite{hemker2003two}, \cite{hemker2004fourier}, and \cite{hemker2004fourierlinear}.  However, convergence rates were in the range of $0.4$ to $0.6$ for low order discretiztions $(p\le 2)$.  In addition, local Fourier analysis is applied to a two level multigrid scheme, and is used as a heuristic to estimate GMG performance.  Further, it is well known that two grid optimality does not always imply V-cycle optimality (see \cite{napov2010does}).
\\
\indent  It should be noted that the HDG class of discretizations is quite large, similar to that of standard DG methods.  In \cite{CockburnGL09}, a unified framework is developed (similar to the work of Arnold et al. in \cite{ArnoldBCM02}) to create a taxonomy of HDG methods.  For instance, one can obtain HDG methods by utilizing one of: Raviart--Thomas DG, Brezzi--Douglas--Marini DG, local DG, or interior penalty DG.  In this work we have no need to distinguish between the various HDG methods, because we employ the LDG family of hybridizable methods (\cite{CockburnDGRS09}).  As such, since the LDG method is a mixed technique, one needs to reformulate the underlying equation \eqref{eq:pde1} as a first order system  by introducing an auxiliary variable ${\bf q}$:
\begin{align}
{\bf q} & = -{\bf K} \nabla u,&&\quad\mbox{in}\quad\Omega,\\
		{\bf K}^{-1} {\bf q} + \nabla u &= 0,&&\quad\mbox{in}\quad\Omega,
		\notag
		\\
		\nabla \cdot   {\bf q}   &= f,&&\quad\mbox{in}\quad\Omega,
		\notag
		\\
		u &= g_D,&&\quad\mbox{on}\quad\partial\Omega.
		\label{eq2}
\end{align}
\indent  We now describe the HDG method.  Let $\mathcal{E}_h$ be a subdivision of $\Omega$, made of quadrilaterals, $K$,  of maximum diameter $h$.  The unit normal vector outward of $K$ is denoted by $\bm n_K$.
The mesh skeleton is denoted by $\Gamma_h$, that is the union of all the edges.  We further decompose the mesh skeleton as $\Gamma_h = \Gamma_h^o\cup\Gamma_h^\partial, $ where $\Gamma_h^\partial$ denotes the set of all edges on the boundary of the domain, and $\Gamma_h^o$ the set of all interior edges.  The broken Sobolev space is represented by $H^1(\mathcal{E}_h)$; it consists of piecewise $H^1$ functions on each mesh element.  We use the following short-hand notation for $L^2$ inner-product on mesh elements and edges:
\begin{align}
		(u,v)_{\mathcal{E}_h} & = 
		\sum_{K \in {\mathcal{E}_h} }\int_K  u v\,dx,
\quad
		\langle u,v\rangle_{\Gamma_h}  = 
		\sum_{K \in {\mathcal{E}_h} }\int_{\partial K}  u v\,ds,
 		&&\forall u,v \in  H^1(\mathcal{E}_h),
\end{align}
\begin{align}
  		\langle {\bm w}\cdot{\bm n}, v \rangle_{\Gamma_h} & = 
		\sum_{K\in\mathcal{E}_h} \int_{\partial K} {\bm w}|_K\cdot{\bm n_K}  v|_K ds,
    		&&\forall {(\bm w},v)\in H^1(\mathcal{E}_h)^2\times H^1(\mathcal{E}_h).
\label{eq:short2}
\end{align}
\noindent  The underlying approximation spaces for the HDG method are as follows:
	\begin{align*}
		{  W}_h & =  \{ w \in  L^2(\Omega):\, w|_K \in \mathbb{Q}_p(K)~~\forall K\in \mathcal{E}_h \}, \quad
{\bm V}_h = W_h\times W_h,
		\\
		{  M}_h & = \{ \mu \in L^2(\Gamma_h):\, \mu |_e \in \mathbb{Q}_p( e )~~\forall e\in \Gamma_h\},
		 \quad {  M}_h^p = {  M}_h,
	\end{align*}
where $\mathbb{Q}_p(K)$ is the standard finite element space for quadrilaterals.  That is, $\mathbb{Q}_p(K)$ is the tensor product of polynomials of degree $p$ on each variable.
The same definition \eqref{eq:short2} applies to $\langle {\bm w}\cdot{\bm n}, \mu\rangle_{\Gamma_h}$ for functions ${\bm w}\in H^1(\mathcal{E}_h)^2$ and $\mu \in M_h$.
\\
\indent  The HDG method seeks an approximation $({\bm q}_h,u_h,\lambda_h)\in  {\bm V}_h\times {  W}_h\times{  M}_h $ of
the exact solution $({\bm q}|_{\Omega}, u|_{\Omega}, u|_{\Gamma_h\setminus\partial\Omega_D})$ such that
	\begin{align}
		( {\bm q}_h,{\bm v})_{\mathcal{E}_h} -(u_h,\nabla \cdot {\bm v})_{\mathcal{E}_h}
+\langle \lambda_h,{\bm v}\cdot {\bm n} \rangle_{\Gamma_h^o} &=
		-\langle P_h g_D, {\bm v}\cdot {\bm n} \rangle_{\partial \Omega},
		\label{eq_h1}
		\\
		-( {\bm q}_h, \nabla w)_{\mathcal{E}_h} 
+ \langle {\bm q}_h \cdot {\bm n}, w\rangle_{\Gamma_h}  + \langle\tau (u_h - \lambda_h),w \rangle_{\Gamma_h^o}
+\langle \tau u_h,w\rangle_{\partial\Omega} &= (f,w)_{\mathcal{E}_h} + \langle \tau P_h g_D,w\rangle_{\partial \Omega},
		\label{eq_h2}
		\\
		\langle {\bm q}_h \cdot {\bm n} ,\mu \rangle_{\Gamma_h} 
+ \langle \tau (u_h-\lambda_h), \mu\rangle_{\Gamma_h^o}
+\langle \tau u_h, \mu\rangle_{\partial\Omega}
&= \langle \tau P_h g_D,\mu\rangle_{\partial\Omega},
		\label{eq_h3}
	\end{align}
for all $({\bm v},w,\mu) \in  {\bm V}_h\times {  W}_h\times{  M}_h $.  The factor $\tau$ is a local stabilization term that is piecewise constant on $\mathcal{E}_h$.  In the above, $P_h g_D$ is the $L^2$--projection of $g_D$, defined by:
\[
\int_e P_h g_D \, \mu = \int_e g_D\, \mu, \quad \forall \mu \in \mathbb{Q}_p(e),\quad \forall e\in\partial\Omega.
\]
\indent  Introducing additional unknowns at first glance does not appear to add much benefit.  However, the HDG method allows us to eliminate the unknowns ${\bm q}_h$ and $u_h$ using equations~\eqref{eq_h1} and~\eqref{eq_h3} in an element by element manner to arrive at a weak formulation in terms of $\lambda_h$ only.  
%
%
Let $\bm Q$, $U$ and $\bm \Lambda$ denote the vectors of degrees of freedom for the numerical
solutions  ${\bm q}_h$, $u_h$ and $\lambda_h$ respectively.   Since the solutions ${\bm q}_h$ and $u_h$ are discontinuous,
we can denote by $\bm Q_K$ and $\bm U_K$ the part of the vectors $\bm Q$, $\bm U$ corresponding to the degrees of freedom located on $K$.  We also denote by $\bm \Lambda_K$ the part of the vector $\bm\Lambda$ that corresponds to the degrees of freedom located on the boundary of $K$.  We can express equations~(\ref{eq_h1}), (\ref{eq_h2}) in matrix form:
\begingroup
\renewcommand*{\arraystretch}{1.5}
 \begin{align}
 {\bm A}_K {\bm Q}_K  -{\bm B}_K^T {\bm U}_K +  {\bm C}_K{\bm \Lambda_K} &= {\bm R}_K,\quad \forall K \in \mathcal{E}_h,
 \label{static_1}
 \\
  {\bm B}_K {\bm Q}_K  +{\bm D}_K {\bm U}_K +  {\bm E}_K {\bm \Lambda}_K &= {\bm F}_K, \quad \forall K \in \mathcal{E}_h.
   \label{static_2}
 \end{align}
\endgroup
\noindent Further, one can obtain the local degrees of freedom:
\begingroup
\renewcommand*{\arraystretch}{1.5}
\begin{equation}
\begin{bmatrix}
\bm Q_K
\\
\bm U_K
\end{bmatrix}
=
\begin{bmatrix}
{\bm A}_K & -{\bm B}_K^T
\\
{\bm B}_K &  {\bm D}_K 
\end{bmatrix}^{-1}
\Bigg(
-
\begin{bmatrix}
{\bm C}_K
\\
{\bm E}_K
\end{bmatrix}
\bm \Lambda_K
+
\begin{bmatrix}
\bm R_K
\\
\bm F_K
\end{bmatrix}
\Bigg)
.
\label{shur0}
\end{equation}
\endgroup
The above inverse is well defined and elimination of the local degrees of freedom for $\bm q$ and $u$  can be done in parallel, independently of one another \cite{CockburnGL09}.  Equation~\eqref{eq_h3} can be written in matrix form.  Fix an interior edge $e\in\Gamma_h$ that is shared by elements
$K_1$ and $K_2$ and denote by ${\bm \Lambda}_{K_1 K_2}$ the set of degrees of freedom for the unknown $\lambda_h$,
that lies on the edge $e$. 
\begin{equation}
{\bm G}_{K_1 e} {\bm Q}_{K_1} 
+{\bm G}_{K_2 e} {\bm Q}_{K_2}
+\tau {\bm H}_{K_1 e}{\bm U}_{K_1}
+\tau {\bm H}_{K_2 e}{\bm U}_{K_2}
-2\tau {\bm M}_e {\bm \Lambda}_{K_1 K_2}
= {\bm 0}.
\label{shur1}
\end{equation}
If the edge $e$ lies on the boundary $\partial\Omega \cap \partial K$ we have
\[
{\bm G}_{K e} {\bm Q}_K
+ \tau {\bm H}_{K e} {\bm U}_K
= {\bm S}_e.
\]

%
%
%
%
\noindent  Now we summarize a general HDG assembly and solve procedure:
\vspace*{3ex}
\begin{enumerate}[(i)]
\setlength\itemsep{0.5em}
\item  Use equation~(\ref{shur0}) to assemble the local problems for $\bm\Lambda_K$.
\item  Assemble the local problems for $\bm\Lambda_K$ into a global system matrix using \eqref{shur1}.
\item  Solve the global system matrix for $\Lambda$.
\item  Using the newly solved for $\Lambda$ in equation~(\ref{shur1}), reconstruct $\bm U$ and $\bm Q$.
\item  Postprocess $U$ and $Q$ to obtain superconvergence.
\end{enumerate}
 
Step (ii) in our framework is technically assembly free, we do not directly assemble and store the global system matrix.  Instead, we exploit the unassembled local problems from equation~(\ref{shur0}) to generate a matrix vector multiplication routine.  The resulting routine utilizes many small but dense matrices, instead of a large sparse matrix.  See \cite{vos2010h} for a survey of different finite element assembly strategies (matrix free or otherwise).  Although the focus of our paper is not assembly, an algorithm was developed in \cite{king2014exploiting} to accelerate the assembly of the HDG global system matrix.

\section{Basis functions}
\label{sec:basis}
We select a tensor product basis for the space $\mathbb{Q}^p(K)$.  To easily facilitate high order approximations, we invoke a nodal basis with nodes that correspond to roots of Jacobi polynomials (Gauss-Legendre-Lobotto (GLL)).  Various operators in the HDG discretization and multigrid method may require evaluation of the nodal basis at points that are not nodal.  Some authors use a modal--nodal transformation to deal with such evaluations (see \cite{karniadakis1999spectral}).  Another technique is to simply use the fast and stable \textit{barycentric interpolation} (see \cite{BerrutT04}).  This allows one to stay in the Lagrange basis and not have to resort to a generalized Vandermonde matrix.  That is, barycentric interpolation allows for the stable evaluation of the Lagrange basis at any point in its domain.  Given $N$ grid points, the setup cost is $\mathcal{O}(N^2)$ to generate the barycentric weights, and a $\mathcal{O}(N)$ cost for each evaluation.  Due do the discontinuous nature of the HDG approximation, evaluations occur element-wise, so the stability of barycentric interpolation is perhaps more important than cost of revaluation.  Since barycentric interpolation allows for arbitrary evaluation, one can select quadrature points that differ from the nodal interpolation points.  Below we introduce the barycentric interpolation formula used in this work.

\begin{figure}[htb!]
\centering
\includegraphics[trim = 10mm 80mm 20mm 85mm, clip, scale = 0.45]{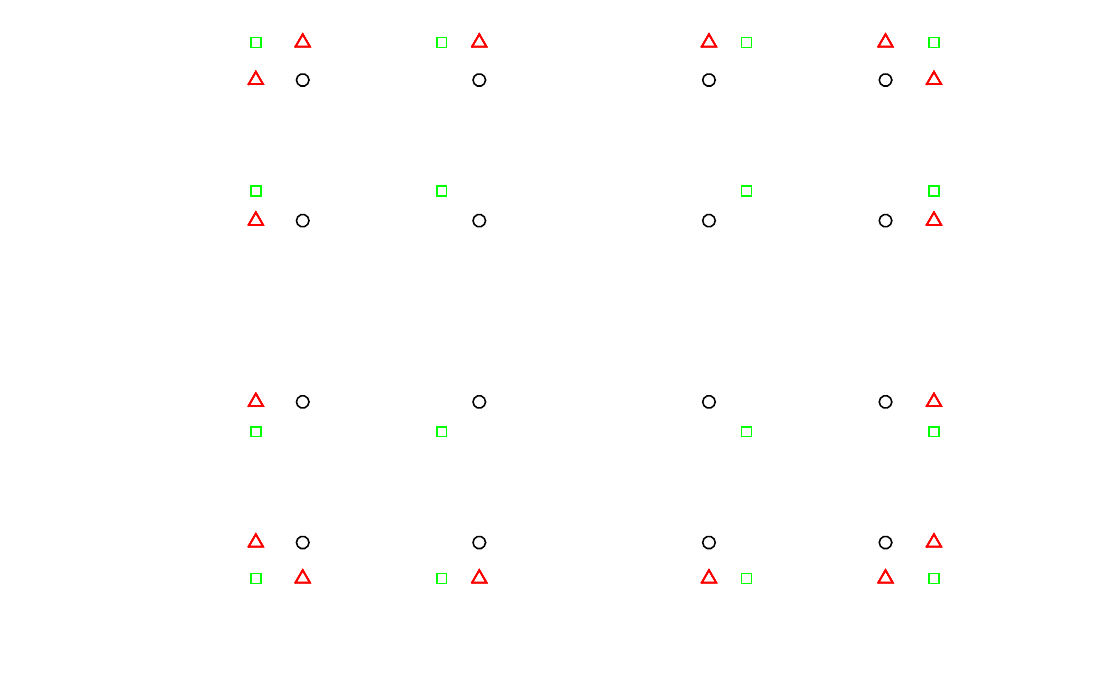}
\caption{GLL nodes (open squares), GL nodes (open circles), GL surface nodes (open triangles).}
\label{sem_grid}
\end{figure}

Let $v$ be a polynomial of degree $n$, that interpolates the scattered data $\{f_0,f_1,\ldots,f_{n}\}$ through the points $\{x_0,x_1,\ldots,x_{n}\}$.  The second form of the barycentric formula is
 
$$
v(x)
=
\frac{ \displaystyle \sum_{j=0}^{n} \frac{W_j}{x-x_j}f_j }{ \displaystyle \sum_{j=0}^{n} \frac{W_j}{x-x_j}  },
$$
 
where $W_j$ are the barycentric weights.  Interpolation points $x_i$ arising from classical orthogonal polynomials have explicit formulas for their barycentric weights (\cite{trefethen2013approximation}, \cite{wang2014explicit}).  As an example, GLL and GL nodes have the following barycentric weights:
 
$$
W_j^{\text{GLL}} = (-1)^j\sqrt{w_j^{\text{GL}}}
~~~
W_j^{\text{GL}} = (-1)^j\sqrt{x_j^{\text{GL}} w_j^{\text{GL}}}
,
$$
 
where $w_j$ and $x_j$ (with appropriate superscripts) are the GL and GLL quadrature nodes and weights.  On the reference element, nodal basis functions are defined using a tensor product of a 1D spectral grid.  In the case of GLL nodes, $\{\xi_i^{(p)}\}_{i=0}^p$ are zeros of a particular family of Jacobi polynomials (see \cite{karniadakis1999spectral}), where $-1\le \xi_i\le 1$.  We adopt the standard that the reference element in is $[-1,1]$ in 1D, $[-1,1]\times [-1,1]$ in 2D, and $[-1,1]\times [-1,1]\times [-1,1]$ in 3D.  See Fig.~\ref{sem_grid} for a sample spectral grid.
 In 2D the tensor product can be written as
  
$$
(\xi_i, \eta_j) \stackrel{def}{=} (\xi_i^{(p)},\eta_j^{(p)}),~~~i,j=0,1,\ldots, p,
$$
 
where $p$ is the polynomial degree associated with $Q^p$.  On each element, we define the basis functions as tensor products:
  
$$
\phi^e_I(\xi,\eta) = \ell_i(\xi)\ell_j(\eta),
$$
 
with $I=i+j(P+1)$ (lexicographic ordering, and $\ell_i$ is the Lagrange polynomial that is nodal at GLL node $i$).  In one dimension, if there are $p+1$ GLL nodes, then there are $p+1$ associated basis functions.  In two dimensions, we take the tensor product of the 1D GLL nodes, which results in $(p+1)^2$ GLL nodes, and $(p+1)^2$ associated basis functions.  In the index $I$, if $i=j=p$, then $I=p+p(p+1)=p^2+2p$, and including the index contribution from $i=j=0$, we have a total of $p^2+2p+1=(p+1)^2$ basis functions in 2D as well.  For brevity, let the lexicographic ordering $I=i+j(p+1) =: (i,j)$.  Let the reference element be given by $\Box = [-1,1] \times [-1,1]$.  Just as in the 1D case, we assume that 
\begin{align*}
u(x,y) \bigg|_{\Omega^e} &\approx \sum_{i=0}^{p}\sum_{j=0}^{p}{u_{ij}^e \phi_I^e( \xi,\eta ) },~~~ ( \xi,\eta )\in \Box,
\\
\phi_I^e(\xi_m,\eta_n) &= \delta_{im}\delta_{jn},~~~i,j,n,m \in \{0,\ldots p \},
\\
(\xi,\eta) &= \text{ GLL quadrature points } \in \Box.
\end{align*}
Using a change of variables to map from the physical element to the reference element we have
\begin{align*}
\iint_\Box u(\xi,\eta) \, d\xi d \eta &= \int_{-1}^1 \bigg(\int_{-1}^1   u(\xi,\eta) \, d\xi \bigg)d \eta
\\
&\approx \sum_{i=0}^P \sum_{j=0}^P w_i w_j u(\xi_i,\eta_j),
\end{align*}
where $w_i,w_j$ are the GLL quadrature weights in the $x$ and $y$ directions.  We want to map an arbitrary element to this element, which requires a change of variables $\xi = \xi(x,y),$ $\eta =\eta(x,y)$.  The Jacobian of this mapping is
$$
\mathcal{J}^e(\xi,\eta)
=
\begin{vmatrix}
 \frac{\partial x}{\partial \xi} & \frac{\partial x}{\partial \eta}
\\
 \frac{\partial y}{\partial \xi}       & \frac{\partial y}{\partial \eta}
\end{vmatrix}
=
\bigg| \frac{\partial x}{\partial \xi}\frac{\partial y}{\partial \eta}- \frac{\partial y}{\partial \xi}\frac{\partial x}{\partial \eta}\bigg|.
$$
Hence, integrating over an arbitrary element $\Omega^e$ ($\cup_{e=1}^E \Omega^e $),
\begin{align*}
\iint_{\Omega^e} u(x,y)\,dx\,dy &= \iint_\Box u^e(\xi,\eta) \mathcal{J}^e(\xi,\eta) \, d\xi \, d\eta
\\
&\approx  \sum_{i=0}^p \sum_{j=0}^p w_i w_j u^e(\xi_i,\eta_j) \mathcal{J}^e_{ij},
\end{align*}
where $ \mathcal{J}^e(\xi_i,\eta_j)= \mathcal{J}^e(\xi,\eta)$.  The above approach assumed that the interpolation points and quadrature points were the same (classic spectral element method).  However, GLL quadrature rule for $p+1$ points is only exact for polynomials of degree $2p-1$.  If higher order quadrature is needed, the GL quadrature rule for $p+1$ points is exact for polynomials of degree $2p+1$.  Then, one can fix the GLL interpolation points, and evaluate the Lagrange basis functions at GL quadrature points if desired.  To do this, we utilize the barycentric formula:
$$
\ell_j(x) = 
\frac{
\displaystyle \frac{W^{\text{GL} }_j}{x-x_j^{\text{GLL} }}
}
{
\displaystyle \sum_{k=0}^{p} \frac{W^{\text{GL}}_k}{x-x_k^{\text{GLL} }}
}.
$$
Barycentric interpolation enables the stable and fast evaluation of the Lagrange polynomial basis $\ell_j$ anywhere in its domain.

\section{HDG discretization}
\label{sec:disc}
 Before proceeding with the results of the HDG GMG method, we verify numerically that the HDG discretization provides the expected optimal $L^2$ convergence rates; $p+1$ for \textit{both} the potential $u_h$ and its flux ${\bm q}_h$.  Moreover, with the use of a local postprocessing filter \cite{CockburnDG08}, we can achieve superconvergence of the potential $u_h$, so that it converges in the $L^2$ norm with the rate $p+2$.  HDG does fall under the umbrella of stabilized DG methods, so the parameter $\tau$ in equations~(\ref{eq_h2}) and~(\ref{eq_h3}) needs to be specified.  The local stability parameter is piecewise constant, defined facet-by-facet.  For the model problem we set $\Omega = [0,1]\times [0,1]$, $\partial\Omega_D=\partial\Omega$ ($\partial\Omega_N=\emptyset$), ${\bf K}(x,y) = \tanh{(x+y)}+1$.  The domain $\Omega$ is partitioned into $N\times N$ squares.  A manufactured solution is used to examine the error: $u(x,y) =    x(x-1)y(y-1)\exp(-x^2-y^2)$, and the corresponding forcing function $f$ is determined from $u$.
\\ 
\indent For the model Poisson problem it turns out that $\tau \equiv 1$ on every facet provides optimal convergence rates.  In Table~(\ref{tab1}), it is apparent that the expected convergence rates are met.  Moreover, the postprocessed potential $u_h^*$ results in a rate of $p+2$.  The postprocessed flux ${\bm q}_h^*$ converges in a rate of $p+1$, but the errors are smaller than that of ${\bm q}_h$. 
 \begin{table}[ht!]
\begin{tabular}{l  b | a | b | a | b | a | b | a | b}
\hline
& &  \multicolumn{2}{c}{$\|u_h-u\|_{L^2(\Omega)}$} & \multicolumn{2}{c}{$\|u_h^*-u\|_{L^2(\Omega)}$} 
&  \multicolumn{2}{c}{$\|q_h-q\|_{L^2(\Omega)}$} & \multicolumn{2}{c}{$\|q_h^*-q\|_{L^2(\Omega)}$} 
\\ 
\cline{3-4} \cline{5-6} \cline{7-8} \cline{9-10}
$p$  &$N$& Error & Rate & Error & Rate&  Error& Rate& Error& Rate
\\\hline
   6 &    2 &  1.05e-07 & -        & 3.02e-09 & -        & 2.39e-07 & -        & 1.12e-07 & -          \\ 
     &    4 &  8.51e-10 & 6.95 & 1.13e-11 & 8.06 & 1.95e-09 & 6.94 & 8.41e-10 & 7.06   \\
     &    8 &  6.97e-12 & 6.93 & 4.45e-14 & 7.99 & 1.61e-11 & 6.92 & 6.65e-12 & 6.98   \\
     &   16 &  5.60e-14 & 6.96 & 7.12e-16 & 5.97 & 1.41e-13 & 6.84 & 6.36e-14 & 6.71   \\
     &   32 &  8.79e-16 & 5.99 & 5.35e-16 & 4.12 & 2.11e-13 & -5.82 & 1.47e-13 & -1.21 \\
\\
   5 &    2 &  7.53e-07 & -        & 3.50e-08 & -        & 1.77e-06 & -        & 9.80e-07 & -        \\
     &    4 &  1.69e-08 & 5.47 & 3.03e-10 & 6.85 & 3.91e-08 & 5.50 & 1.81e-08 & 5.76\\
     &    8 &  2.90e-10 & 5.87 & 2.43e-12 & 6.96 & 6.72e-10 & 5.86 & 2.95e-10 & 5.94\\
     &   16 &  4.72e-12 & 5.94 & 1.92e-14 & 6.98 & 1.10e-11 & 5.94 & 4.69e-12 & 5.97\\
     &   32 &  7.52e-14 & 5.97 & 3.68e-16 & 5.70 & 2.05e-13 & 5.74 & 1.01e-13 & 5.53\\
\\
   4 &    2 &  1.46e-05 & -        & 7.04e-07 & -        & 3.32e-05 & -        & 1.79e-05 & -        \\
     &    4 &  4.92e-07 & 4.89 & 1.09e-08 & 6.02 & 1.13e-06 & 4.88 & 5.49e-07 & 5.03\\
     &    8 &  1.65e-08 & 4.90 & 1.73e-10 & 5.97 & 3.80e-08 & 4.89 & 1.76e-08 & 4.96\\
     &   16 &  5.37e-10 & 4.94 & 2.75e-12 & 5.98 & 1.24e-09 & 4.94 & 5.62e-10 & 4.97\\
     &   32 &  1.71e-11 & 4.97 & 4.33e-14 & 5.99 & 3.98e-11 & 4.97 & 1.78e-11 & 4.98\\ 
\\     
   3 &    2 &  8.61e-05 & -        & 7.58e-06 & -        & 2.05e-04 & -        & 1.35e-04 & -        \\
     &    4 &  7.78e-06 & 3.47 & 2.63e-07 & 4.85 & 1.80e-05 & 3.52 & 9.86e-06 & 3.78\\
     &    8 &  5.49e-07 & 3.82 & 8.63e-09 & 4.93 & 1.27e-06 & 3.82 & 6.59e-07 & 3.90\\
     &   16 &  3.63e-08 & 3.92 & 2.77e-10 & 4.96 & 8.43e-08 & 3.91 & 4.28e-08 & 3.95\\
     &   32 &  2.33e-09 & 3.96 & 8.80e-12 & 4.98 & 5.43e-09 & 3.96 & 2.73e-09 & 3.97\\
\hline
 \end{tabular}
\caption{Errors and convergence rates of the HDG scheme, on a Cartesian mesh of $N\times N$ elements.}
\label{tab1}
\end{table}

\section{Multigrid algorithm}
\label{sec:mg_description}
Multigrid methods work by eliminating a wide range of frequencies from the estimated error between the exact solution, and the multigrid iterate.  Any multigrid method has three key components: relaxation, grid transfers, and coarse grid operators (course grid correction).  In the following subsections our selection of these components is presented.  The geometric multigrid method leverages information about the underlying geometry, discretization, and PDE.  This results in an optimal solver (\cite{trottenberg}, \cite{Brandt2}) that is coupled to the discretization; only $\mathcal{O}(N)$ floating point operations are required to reduce the error to discretization level accuracy.  Moreover, GMG offers a uniform convergence rate.  That is, the rate of convergence of GMG is not dependent on the mesh size.
\\
\indent In direct contrast to GMG, algebraic (black box) solvers rely entirely on the discretization matrix.  The trade off with algebraic solvers (e.g. Krylov methods) is that they are decoupled from the discretization, but their convergence is very sensitive to the condition number of the discretization matrix.  Refinement in $h$ or $p$ increases said condition number, and in turn the algebraic solver requires even more iterations to reach a desired tolerance.  In this context, very efficient preconditioners are needed to assist the algebraic solvers.  Multigrid methods are also very popular as preconditioners (\cite{gopalakrishnan2003multilevel}).
\\
\indent  The finite element method brings with it two ways to potentially increase accuracy: $h$ refinement and $p$ refinement.  Since GMG is coupled to the discretization, it can leverage $h$ refinement and/or $p$ refinement (refinement in $p$ means increasing/decreasing the polynomial order).  This implies that two different classes of multigrid operators are needed.  One for $h$ refinement, the other for $p$ refinement.  For simplicity, we do not consider $hp$--GMG (refining $h$ and $p$ simultaneously).  In $h$--GMG the polynomial order $p$ is fixed, and the hierarchy of grid spaces are determined by increasing/decreasing the mesh size.  In $p$--GMG the mesh size $h$ is fixed, and the hierarchy of grid spaces are determined by their increasing/decreasing the polynomial order.  The combination of these three refinement strategies offers interesting and potentially efficient ways to deal with high order and very high order approximations (\cite{Ant}, \cite{Fidkowski}). 
\\
\indent We use a nested multigrid strategy.  The multigrid method starts with $p$-multigrid, which fixes the mesh and the fine and coarse grids are formed by increasing/decreasing the polynomial order.  Once the polynomial order reaches $p=0$, which corresponds to a face-centered finite volume scheme, we continue with the $h$-multigrid scheme as suggested in \cite{gopalakrishnan2009convergent}, \cite{cockburn2013multigrid}, \cite{tan2009iterative}.  This technique exhibits $h$-independent convergence.  It projects the HDG $p=0$ approximation into an appropriate space of continuous functions, so that standard $h$-multigrid for continuous Galkerin can be used for the coarser levels.  A diagram of our general nested multigrid strategy is given in Fig.~\ref{fig_mg_description}.
\\  
\indent  We use rediscretization for $p$-multigrid to form the hierarchical grids; this way we do not store the sparse linear systems for the polynomial orders $p>0$.  The corresponding grid transfer operators are also implemented in a matrix free manner, as the canonical grid transfer operators are discontinuous across element interfaces, they are completely data parallel.
\\
\indent  For $p=0$, and subsequent $h$-multigrid levels, we store the linear systems in compressed row storage format.  On the coarsest level (matrix of about dimension 2000) of the $h$-multigrid scheme, we employ the MKL optimized PARDISO direct solver \cite{schenk2004solving}.  Other options are available for the coarse grid solve, for instance replacing the direct solve with a number of relaxation steps, or continuing the multigrid hierarchy with an algebraic multigrid method.
\\
\indent  For problems where geometric multigrid is not suitable, or generating a meaningful geometric hierarchy is not possible, one can replace $h$-multigrid with algebraic multigrid or an appropriate coarse grid solver.  In these cases, $p$-multigrid is still applicable as the mesh is fixed.

\begin{figure}[ht!]
\centering
\begin{tikzpicture}[ scale = 0.8, rotate = 180 ]
\coordinate (A) at (-4.5,0) {};
\coordinate (B) at ( 4.5,0) {};
\coordinate (C) at (0,8) {};
\coordinate (C1) at (1.1,6) {};
\coordinate (C2) at (-1.1,6) {};

\coordinate (O1) at (0,1.3) {};
\coordinate (O2) at (0,3.3) {};
\coordinate (O3) at (0,5.3) {};
\path[name path=AC,draw=none] (A) -- (C);
\path[name path=BC,draw=none] (B) -- (C);
\filldraw[draw=black, ultra thick,fill=blue!10] (A) -- (B) -- (C1) -- (C2) --cycle ;

\filldraw[black!20, draw=black, thick] (O1) node[black, above] {$p$-multigrid};
\filldraw[black!20, draw=black, thick] (O2) node[black, above] {$h$-multigrid};
\filldraw[black!20, draw=black, thick] (O3) node[black, above] {Direct solve};

\draw[draw=red, very thick, dashed] (2.2,4) -- (-2.2,4);
\draw[draw=red, very thick, dashed] (3.2,2) -- (-3.2,2);
\end{tikzpicture}
\caption{Diagram of our $hp$-multigrid algorithm.}
\label{fig_mg_description}
\end{figure}

\subsubsection{$p$-multigrid grid transfer}
Finite element methods pair well with multilevel methods, due to the flexibility in chose of underlying spaces.  Consider the nested discontinuous finite element spaces $\Omega^H\subset \Omega^h$.  The mapping from a coarse space ($\Omega^H$) to a fine space ($\Omega^h$) is called prolongation, or interpolation.  There are many ways to build a prolongation mapping, however, in the context of finite element methods there is a canonical choice. 

The canonical prolongation operator uses the fact that $\Omega^H\subset \Omega^h$.  That is, any coarse grid function can be expanded as a linear combination of fine grid functions.  Thus, we take the prolongation operator to be the natural embedding from  $\Omega^H$ into $\Omega^h$.  As such, $P:\Omega^H \to \Omega^h,$ 
\begin{equation}
 P(u_H)=u_H\quad \forall u_H \in \Omega^H.
\label{eq_mg_pro}
\end{equation}
 The mapping from a fine space to a coarse space is called restriction, or sub-sampling.  It is not as trivial to select a restriction operator.  However, since finite element methods are based on a weak formulation, a canonical restriction is defined by a $L^2$ projection.  The restriction is given by $R: (\Omega^h)'\to (\Omega^H)' $, 
\begin{equation}
(R v_h,u_H) = (v_h,u_H),~~~\forall v_h \in (\Omega^h)',~~~\forall u_H\in \Omega^H .
\label{eq_mg_res}
\end{equation}
Notice that $(R v_h,u_H) =(  v_h,R^*u_H)  = (v_h,Pu_H)$.  One can infer from this that $R^*=P$.  Selecting the restriction operator in this way is deliberate; it is done to ensure that the multigrid operator is symmetric.  This allows the multigrid method to be a suitable preconditioner for symmetric iterative methods.
\\
\indent  The grid transfer operators are defined on the reference element, and for discontinuous Galkerin methods the canonical prolongation and restriction operators are data independent.  For $p$-multigrid we apply the grid transfer operators in a matrix free manner.  These grid transfer operators are defined in equations~\eqref{eq_mg_pro} and~\eqref{eq_mg_res}, where it is understood that for $p\ge 0$, the coarse space is $\Omega^H={  M}_h^{p}$, and the fine space is $\Omega^h= {  M}_h^{p+1}$.
\\
\indent  We note that for the HDG method the grid transfer operators are defined on a lower dimension.  This is due to static condensation reducing the unknowns (dimension $d$) to the mesh skeleton (dimension $d-1$).  For instance, classical DG methods in 2D would have a prolongation operator $P:Q^p\to Q^{p+1}$ of size $(p+2)^2\times (p+1)^2$.  For HDG in 2D, the prolongation operator would be $P$ of size $(p+2)\times (p+1)$, as static condensation reduces the unknowns to 1D.  Nodal basis functions give rise to dense grid transfer operators, whereas modal basis functions have sparse binary grid transfer operators.
\subsubsection{$h$-multigrid grid transfer}
The $h$-multigrid scheme we use is taken from \cite{gopalakrishnan2009convergent}, \cite{cockburn2013multigrid}, \cite{tan2009iterative}.  The main idea is to project the $p=0$ HDG solution into the space of continuous functions, then continue with a standard geometric multigrid for piecewise linear continuous finite elements.  It is also possible to project the $p=1$ HDG solution to the piecewise linear finite element space.  However, in~\cite{tan09} and \cite{Cockburn17102013}, it is found that this leads to larger iteration counts and poor convergence rates compared to the $p=0$ projection.  Our numerical experiments also agree with their findings.
\\ \noindent
Let $\mathcal{E}_{k}$ denote a uniform mesh of quadrilateral elements at level $k$.  The coarsest mesh is denoted by $k=0$.  Geometric refinements of the coarsest mesh are then made, for a total of $J$ meshes.  The mesh for the original discretization (finest mesh) is denoted by $k=J-1$.  For simplicity we use uniform meshes so that a geometric hierarchy can be used.  Set 
$$
M_k = \{v: \Omega \to \mathbb{R}: v \text{ is continuous } v |_{\partial\Omega}=0,~v|_K \in Q^1(K), ~\forall K \in \mathcal{E}_{k+1}\},
$$
for $0\le k \le J-1$, and $M_J = M_h^0$.  Note that the spaces have a non-nested property $M_0\subset M_1\subset \ldots \subset M_{J-1} \not\subset M_J.$  As such, a specialized prolongation operator has to be used \cite{gopalakrishnan2009convergent}.  It is defined as $P_k: M_{k-1} \to M_k$ ($2\le k \le J$),
\begin{equation}
P_k v
= 
\begin{cases}
v, &\mbox{if } k<J ,
\\
\Pi_{M_J} (v|_{\Gamma_h}), &\mbox{if } k=J,
\end{cases}
\label{eq_h_mg_pro}
\end{equation}
where $\Pi_{M_J}$ is the $L^2$-orthogonal projection onto $M_J$.  The corresponding restriction operator $R_{k-1}: M_k \to M_{k-1}$ follows from equation~\eqref{eq_mg_pro} ($1\le k \le J-1$), 
\begin{equation}
(R_{k-1} w,u)_{k-1} = (w, P_{k}u)_k,\quad \forall u,w \in M_{k-1}.
\label{eq_h_mg_res}
\end{equation}

 \subsubsection{Coarse grid operator}
The multigrid method utilizes a hierarchy of grids.  As such, there is a need for a hierarchy of discretizations.  In other words, on a coarse grid we have a reduced model of our original discretization.  This reduced model needs an appropriate discretization (coarse grid operator) for the coarse grid.  Two of the most common coarse grid operators are subspace inheritance and subspace non-inheritance.  

\indent Subspace inheritance defines coarse grid operators by ${\bf A}^H = {\bf R} {\bf A}^h {\bf P}$, where ${\bf A}^H$ is the coarse grid operator, ${\bf A}^h$ is the fine grid operators, ${\bf R}$ and ${\bf P}$ are the restriction and prolongation operators, respectively.  On each coarse grid information from the previous level is being inherited.  This technique has the disadvantage that continuity may be implicitly enforced for sufficiently coarse grids.  Subspace non-inheritance simply uses rediscretization, on each coarse grid $\Omega^H$ the discretization operator is built for the dimension associated with $\Omega^H$.

\indent In the $p$-multigrid setting, the coarse grid problems are formed by discretizating our problem for smaller polynomial degrees.  With $p$-multigrid there is a subtle choice that is to be made for its refinement strategy.  The coarse spaces for $p$-multigrid are obtained by decreasing the polynomial order $p$.  How one should decrease the polynomial order is not well agreed upon. A couple of coarsening ratios have been proposed in the context of continuous $hp$-multigrid (see \cite{Mitchell}): 
\begin{itemize}
\item geometric: assumes $p$ is a power of two, $p,p/2,p/4,\ldots,2,1,2,4,\ldots$,
\item odd: assumes $p$ is a an odd natural number,
\item even: assumes $p$ is a an even natural number,
\item full arithmetic: $p,p-1,p-2,\ldots,2,1,0,1,2,\ldots$.
\end{itemize}
The most flexible coarsening ratio is full arithmetic, which allows for arbitrary polynomial degrees. 
\\
\indent  To define the $p$-multigrid hierarchy we use full arithmetic coarsening.  We selected full arithmetic coarsening since we are using multigrid as a solver, not as a preconditioner.  In this case, there are more opportunities to eliminate high frequency modes.  For high order DG methods, traditional parallelizable relaxation schemes do not properly eliminate high frequency modes \cite{johannsen2005multigrid,johannsen2005symmetric,gopalakrishnan2003multilevel,Antonietti2016,brenner2005convergence,brenner2009multigrid}.  Full arithmetic coarsening was found to be as effective as more rapid coarsening for problems with moderate polynomial degrees \cite{Mitchell,Fidkowski}. 
\\
\indent  For $p > 0$, we employ a matrix free procedure.  In the $h$-multigrid setting, the coarse grid problems are formed by subspace inheritance.  For low order problems, there is little benefit for matrix free operations, so we store the $p=0$ matrix and any subsequent matrices that arise from the $h$-multigrid algorithm.  This amounts to the Galerkin triple product ${\bf A}_{k-1}= {\bf R}_{k-1} {\bf A}_k {\bf P}_k$ (for $1\le k \le J-1$), where ${\bf R}_{k-1}$ and ${\bf P}_k$ are as defined in equations~\eqref{eq_h_mg_res} and \eqref{eq_h_mg_pro}, respectively.   

\subsubsection{Relaxation}
Relaxation (also called smoothing) plays an integral part in multigrid techniques.  High frequencies are damped using a relaxation method, leaving low frequencies to be resolved by coarse grid correction.  In some sense, multigrid can be thought of as a way to accelerate the convergence of a relaxation method; typically the relaxation methods used in multigrid have poor (or no) convergence properties.  The main purpose of a relaxation method is to smooth the error, not necessarily make it small.  Standard choices for relaxation are pointwise/block stationary iterative methods (Jacobi, damped Jacobi, Gauss-Seidel, Successive over relaxation, etc.), and polynomial type preconditioners.\\
\indent  A general overview of iterative methods can be found in~\cite{saad}.  For lower order discretizations, (in the context of this work, third order or less) standard relaxation methods tend to work well.  However, for higher order discretizations, loss of uniform convergence is experienced.  A larger number of relaxation steps, or more powerful smoothers are required to retain uniform convergence in the high order case.
\\
\indent  In a parallel setting, polynomial smoothers where shown in~\cite{Adams2003593} to perform better in terms of time-to-solution.  With the advent of many-core devices, some authors have examined parallel asynchronous iterations (see \cite{Anzt2012}, \cite{Anzt2013}).  Sparse approximate inverses (SPAI) have exhibited promising results when used as preconditioners and smoothers in multigrid (see \cite{Benson}, \cite{Broker}).  Moreover, the construction of SPAI operators are inherently parallel, and only require a matrix vector product to apply. 
\\
\indent  Relaxation for high order discontinuous Galkerin methods typically require specialized treatment.  This was found to be the case for HDG in our experiments, as well as in \cite{tan2009iterative} and \cite{yakovlev2016cg}.  In particular, block-type methods work better than pointwise or polynomial smoothers.  In the context of our work, we have found that SPAI and additive Schwarz domain decomposition methods work very well as a smoother.  In Section~\ref{sec:gmg} we present multigrid convergence results for an SPAI preconditioner.  The additive Schwarz smoother has similar convergence properties depending on the overlap size.  These more expensive smoothers are used in our work as they able to retain robust multigrid convergence.  Section~\ref{sec:gmg} discusses this matter further.

\subsubsection{Sparse approximate inverse relaxation}
Sparse approximate inverses have exhibited promising results when used as preconditioners and smoothers in multigrid (see \cite{Benson}, \cite{Broker}).  Moreover, the construction of SPAI operators are inherently parallel, and only require a matrix vector product to apply.  SPAI seeks an approximation ${\bm M} $ of ${\bm A}^{-1}$ such that
 
$$
\min_{ {\bm M} \in \mathcal{S}}\| {\bm I} - {\bm M}{\bm A} \|_F,
$$
 
where $ \mathcal{S}$ is a set of sparse matrices.  By using the Frobenius norm, a great deal of concurrency is exposed:
 
$$
\| {\bm I} - {\bm M}{\bm A} \|_F^2
 =
 \sum_{i=1}^{n} \| \vec{e}_i^T - {\bm M}_i^T{\bm A} \|_2^2.
$$
 
Each of the $n$ least squares problems are independent of one another.  Storage and computational savings can be obtained by leveraging the sparsity of $\bm A$ and $\bm M$.  In addition, since the HDG method gives rise to a symmetric operator, the factorized sparse approximate inverse (FSAI) can be utilized to obtain further improvements.  The use of SPAI as a smoother will require its explicit storage, but the application is that of a matrix vector multiply:
 
$$
\vec{x}_{\text{SPAI}}^{(k+1)} = \vec{x}_{\text{SPAI}}^{(k)} + {\bm M}( \vec{f} - {\bf A}\vec{x}_{\text{SPAI}}^{(k)} ).
$$
 
\subsubsection{Additive Schwarz domain decomposition relaxation}
Increasing the polynomial order in finite element methods causes the condition number of the discretization operator to grow.  This stretches the spectrum of the discretization operator, and traditional smoothers can no longer damp a sufficiently wide range of high frequencies.  To combat this, one can increase the number of smoothing steps, or seek a more powerful smoother.  In some cases, domain decomposition methods can provide a very powerful smoother for multigrid techniques.  In \cite{Lottes2005}, \cite{Stiller15}, and \cite{Stiller16}, variations of the Schwarz domain decomposition method is used as a multigrid smoother for high order finite element discretizations.
\\
\indent  The standard additive Schwarz method takes the form
 
$$
{\bm M} = \sum_{e=1}^{|\Gamma_h|} {\bm R}_e^T {\bm A}_e^{-1} {\bm R}_e,
$$
 
where ${\bm R}_e$ is a binary restriction operator that transfers global datum to local datum.  The operator ${\bm A}_e$ is the discretization operator restricted on the facet $e$.  Notice that this is a single level Schwarz method, and as such, it is not as effective as their two level counter parts.  Weighting the single level Schwarz method has been demonstrated to produce effective smoothers (see \cite{Lottes2005},\cite{Stiller16}, \cite{Loisel2008}).  The weighting appears to be successful, but is also somewhat experimental.  The inverse of a diagonal counting matrix is used in \cite{Lottes2005}, and a type of bump function is used in \cite{Stiller16} (both of these weightings are used without a convergence theory).  For a weighting matrix $\bm W$, the additive Schwarz smoother is applied as follows:
 
$$
\vec{x}_{\text{ASM}}^{(k+1)} = \vec{x}_{\text{ASM}}^{(k)} + {\bm W}{\bm M}( \vec{f} - {\bf A}\vec{x}_{\text{ASM}}^{(k)} ).
$$
 
On each level of the $p$-multigrid hierarchy we store smoother $\bf M$ as dense blocks, and apply them in a unassembled manner.  The size of these blocks depends on the polynomial order and overlap size.
\subsubsection{Standard multigrid V-cycle}
Our nested multigrid V-cycle is presented in Algorithms~\ref{algo_1} and~\ref{algo_2}.  This is often referred to as the $V$-cycle.  We use the following notation: $S^\nu(\cdot )$ denotes $\nu$ steps of a relaxation method.  In the $p$-multigrid cycle the Smoother $S^\nu(\cdot )$ is either SPAI or the additive Schwarz relaxation.  In Section~\ref{sec:gmg} we use SPAI, and in Section~\ref{sec:matvec} we use the additive Schwarz relaxation.  For $h$-multigrid we use the Chebyshev smoother \cite{saad} as the problem is low order.
\begin{center}
 \begin{algorithm}[H]
\begin{algorithmic}[1]
\Statex ${\bf A}^h$, ${\bf R}$, ${\bf P}$ are stored in sparse matrix format
\Statex ${\bf R}$ and $\bf P$ are as defined in equations~\eqref{eq_h_mg_res} and~\eqref{eq_h_mg_pro}
\Statex The Smoother $S$ is a polynomial relaxation (Chebyshev).
\State $v^h \leftarrow S^{\nu_1}({\bf A}^h,r^h,v^h)$
\If {On coarsest level}
\State  $v^H \leftarrow ({\bf A}^{H})^{-1} v^H$ (bottom level solver, e.g. direct solve or AMG)
\Else
\State $r^H \leftarrow {\bf R} ( r^h - {\bf A}^h v^h)$
\State $v^H \leftarrow 0$
\State $v^H \leftarrow hMG(v^H,r^H)$
\EndIf
\State $v^h \leftarrow v^h + {\bf P} v^H$
\State $v^h \leftarrow S^{\nu_2}({\bf A}^h,r^h,v^h)$
\end{algorithmic}
\caption{$v^h\leftarrow hMG(v^h,r^h)$}
\label{algo_1}
\end{algorithm}
\end{center}

\begin{center}
 \begin{algorithm}[H]
\begin{algorithmic}[1]
\Statex ${\bf A}^h$, ${\bf R}$, ${\bf P}$ are applied in a matrix free manner
\Statex ${\bf R}$ and $\bf P$ are as defined in equations~\eqref{eq_mg_res} and~\eqref{eq_mg_pro}
\Statex The Smoother $S$ is either SPAI or the additive Schwarz relaxation.
\State $v^h \leftarrow S^{\nu_1}({\bf A}^h,r^h,v^h)$
\If {On coarsest level ($p=0$))}
\State $v^H \leftarrow hMG(v^H,r^H)$ (bottom level solver, is $h$-multigrid)
\Else
\State $r^H \leftarrow {\bf R} ( r^h - {\bf A}^h v^h)$
\State $v^H \leftarrow 0$
\State $v^H \leftarrow pMG(v^H,r^H)$
\EndIf
\State $v^h \leftarrow v^h + {\bf P} v^H$
\State $v^h \leftarrow S^{\nu_2}({\bf A}^h,r^h,v^h)$
\end{algorithmic}
\caption{$v^h\leftarrow pMG(v^h,r^h)$}
\label{algo_2}
\end{algorithm}
\end{center}

\section{Multigrid convergence}
\label{sec:gmg}
\begin{figure}[ht!]
\centering
    \captionsetup{justification=centering}
    \subfloat[V-Cycle.]{
        \includegraphics[trim = 30mm 85mm 45mm 90mm, clip, scale = 0.45]{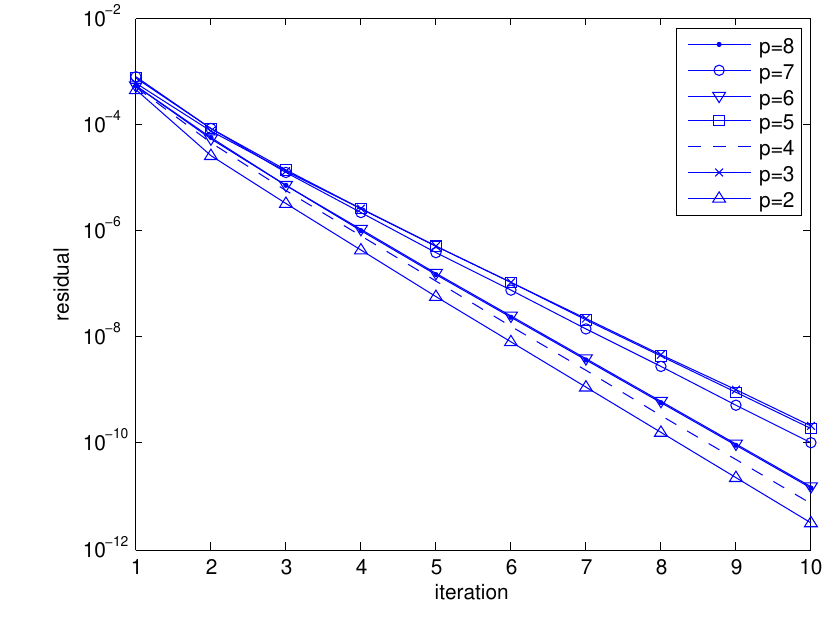}
        \label{fig_gmg1}
    }
    \subfloat[Convergence rate (V-Cycle).]{
        \includegraphics[trim = 30mm 85mm 45mm 90mm, clip, scale = 0.45]{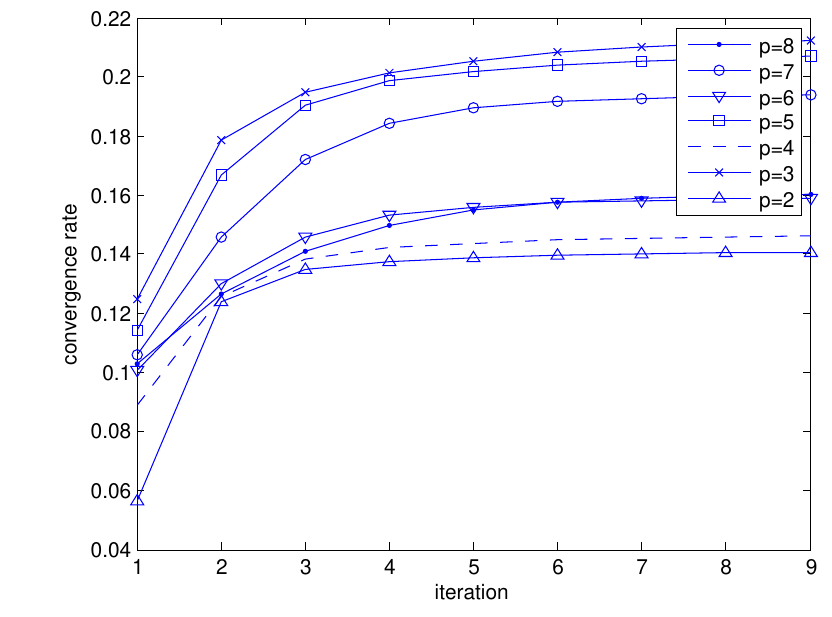}
        \label{fig_gmg2}
    }   
    \\
    \subfloat[W-Cycle.]{
        \includegraphics[trim = 30mm 85mm 45mm 90mm, clip, scale = 0.45]{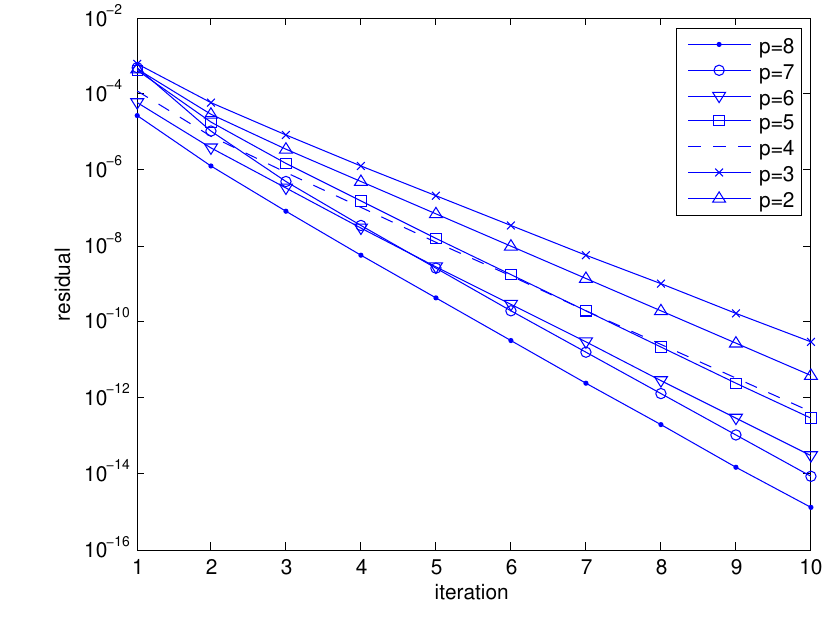}
        \label{fig_gmg3}
    } 
    \subfloat[Convergence rate (W-Cycle).]{
        \includegraphics[trim = 30mm 85mm 45mm 90mm, clip, scale = 0.45]{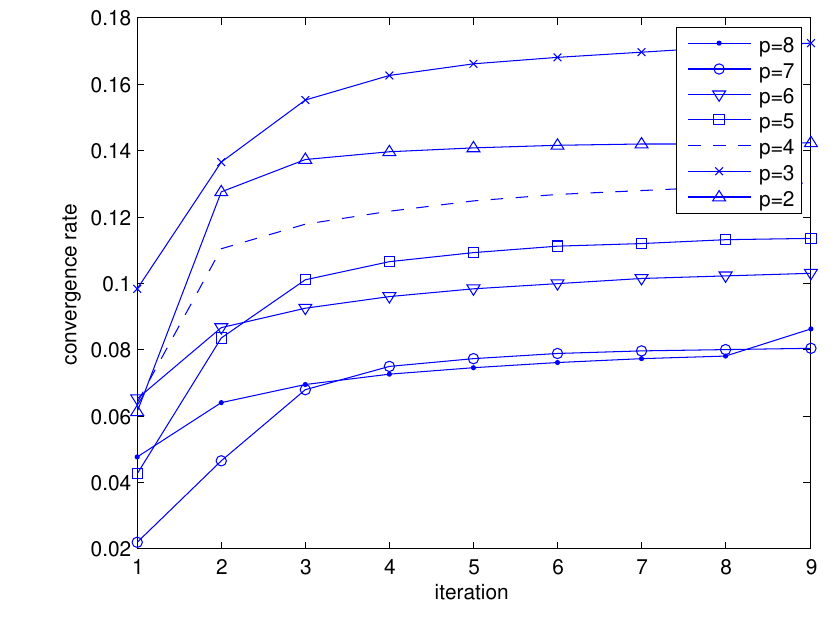}
        \label{fig_gmg4}
    }
\captionsetup{justification=justified}
\caption{GMG for HDG (SPAI-1 smoother).}
\label{fig_gmg5}
\end{figure}

\begin{figure}[ht!]
\centering
    \captionsetup{justification=centering}
    \subfloat[V-Cycle.]{
        \includegraphics[trim = 30mm 85mm 45mm 90mm, clip, scale = 0.45]{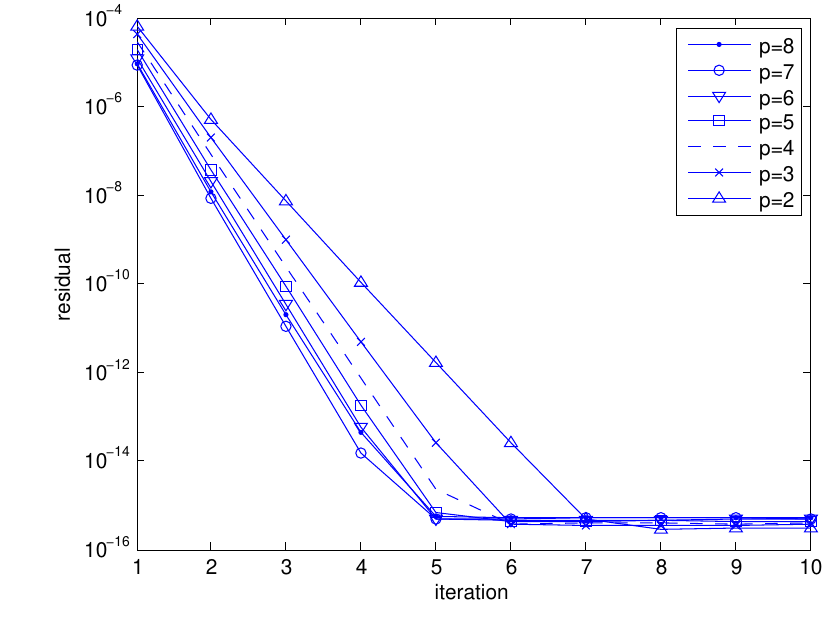}
        \label{fig_fmg1}
    }
    \subfloat[Convergence rate.]{
        \includegraphics[trim = 30mm 85mm 45mm 90mm, clip, scale = 0.45]{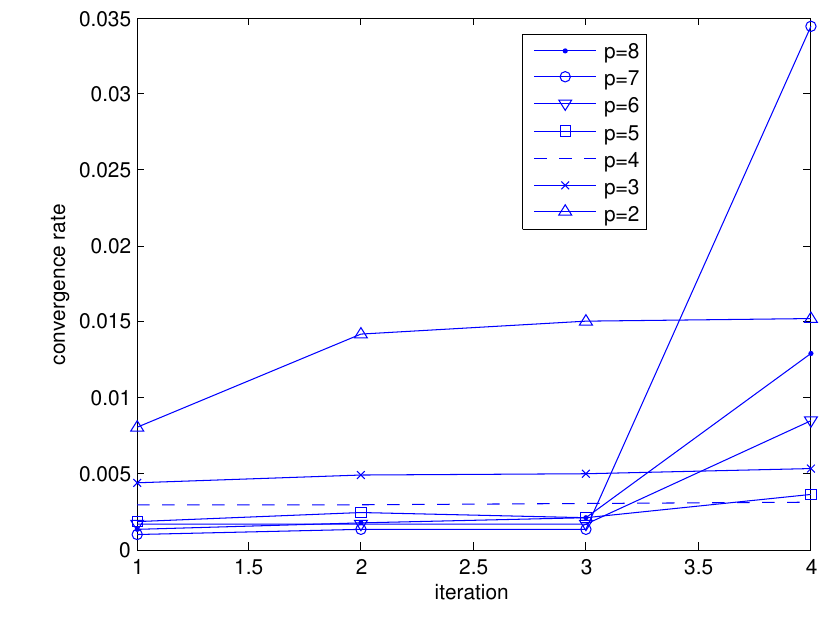}
        \label{fig_fmg2}
    }   
    \\
    \subfloat[Convergence rate.]{
        \includegraphics[trim = 30mm 85mm 45mm 90mm, clip, scale = 0.45]{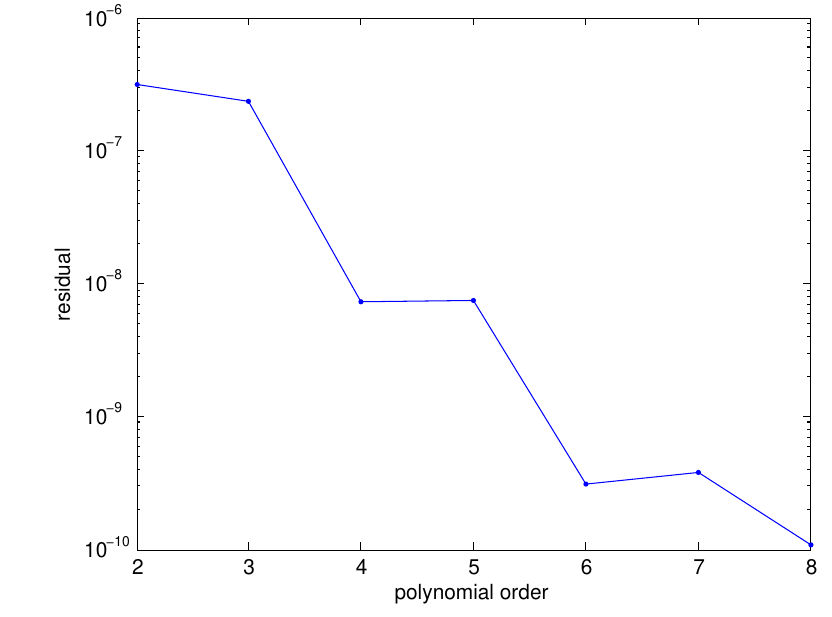}
        \label{fig_fmg3}
    } 
    \subfloat[Operator complexity.]{
        \includegraphics[trim = 30mm 85mm 45mm 90mm, clip, scale = 0.45]{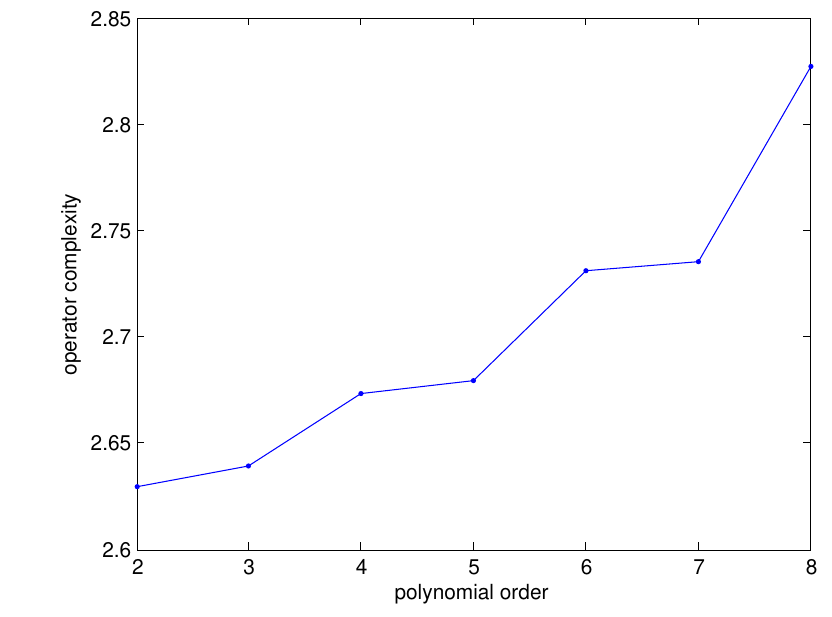}
        \label{fig_fmg4}
    }
\captionsetup{justification=justified}
\caption{FMG for HDG (aggressive FSAI).}
\label{Fig:fig_fmg5}
\end{figure}

In \cite{Ant}, uniform convergence in $h$ and $p$ was established for a discontinuous Galerkin $hp$-multigrid algorithm.  However, in order to guarantee $p$ (or $h$) independent convergence, the number of relaxation steps per level grows quadratically in $p$.  Numerical experiments suggest that for our HDG $p$-multigrid, the relaxation steps per level need to grow linearly (in $p$) to guarantee $p$ independent convergence.  In \cite{cockburn2013multigrid}, $h$ independent convergence was proven for the HDG $h$-multigrid algorithm we use. 
\\
\indent  Fig.~\ref{fig_gmg1} and Fig.~\ref{fig_gmg2} display the results for $2\le p \le 8$ and fine mesh with $(2^5)^2$ elements.  We employ a FSAI smoother on each level, with $\nu_1=\nu_2 = 2$ pre and post smoothing steps.  The FSAI smoother is constructed so that it results in an approximate inverse with a operator complexity of unity (nnz[$\bm M$]$=$nnz[$\bm A$]).  Subspace non-inheritance is used to generate coarse grid operators.  For the $p$-GMG phase of the GMG method, we use full arithmetic coarsening.  To measure the convergence rate, we keep track of the two norm of the fine grid residual between successive iterations:
 $$
 \rho_k = \frac{ \| (r^h)^{(k)} \|_2 }{  \| (r^h)^{(k-1)} \|_2 }.
 $$
 We can see that the results are quite good, even for a modest FSAI smoother.  All of the convergence rates are under 0.22, and tend to cluster in the range 0.13-0.15 for even polynomial degrees.  This observation perhaps indicates that even-$p$ coarsening or geometric-$p$ will not only be more efficient (more rapid coarsening), but also more accurate. 
 \\
 \indent   Fig.~\ref{fig_gmg5} seems to suggest that our method is not exhibiting a $p$ independent convergence rate.  However, we recall that according to \cite{Ant}, uniform convergence in $p$ guaranteed only if the number of relaxation steps per level grows quadratically in $p$ (we set $\nu_1=\nu_2=2$).  For multigrid methods there is a subtle observation to be made: uniform convergence is an excellent property; but if the uniform convergence rate is near one, this is not an efficient solver.  Moreover, if the number of smoothing steps is sufficiently large, a great deal of efficiency is lost.  For reasonable smoothing steps, in \cite{brenner2005convergence} and \cite{Ant} the authors demonstrate convergence rates of 0.85 (or worse).  In order to drop the convergence rates down to $0.5$, upwards of 20 smoothing steps every level needs to be taken.
 \\
\indent  What is needed for efficiency is as follows: uniform convergence, a fast convergence rate (problem dependent), and a small number of uniform smoothing steps.  Satisfying these three properties in practice for high order discontinuous Galkerin methods is not trivial.  Some potential avenues to explore are stronger (but more expensive) smoothers, different coarse grid transfer operators, different multigrid schedules, as well as different grid hierarchies.  This direction of inquiry is beyond the scope of our paper.  In Fig.~\ref{Fig:fig_fmg5} we do show that a more expensive smoother essentially yields uniform convergence with a fast convergence rate and a fixed number of smoothing steps.  
 \\
\indent  For the next set of experiments, we allow the FSAI operator complexity to grow ((nnz[$\bm M$]$\,\approx\,2.5\,$nnz[$\bm A$]).  In Figs.~\ref{fig_fmg1}, \ref{fig_fmg2}, \ref{fig_fmg3}, and \ref{fig_fmg4}, display the results of this change.  For the $V$-cycle, the residual is reduced to machine precision after only 5 to 7 iterations (Fig.~\ref{fig_fmg1}).  The aggressive smoothing also yields stellar convergence rates, below 0.015, as can be ascertained from Fig.~\ref{fig_fmg2}.
\\
\indent  The performance of the V-cycle is better than ideal, so we can easily extend our HDG GMG method to leverage the \textit{full multigrid cycle} (FMG).  FMG is well known to be \textit{the} optimal multigrid schedule for linear problems; with a single FMG iteration, the residual is reduced to discretization level error.  Moreover, it only requires $\mathcal{O}(N)$ floating point operations to achieve this accuracy (\cite{Brandt1}, \cite{Brandt2}, \cite{trottenberg}).  Indeed, Fig.~\ref{fig_fmg3} numerically verifies that a \textit{single} FMG iteration is enough to reach discretization level error.  There is something that particularly interesting about Fig.~\ref{fig_fmg3} - the FMG iteration performs better for higher orders $6\le p$.  This second experiment of course comes at a price: the FSAI smoothers on each level have an operator complexity in the rage of 2.65 to 2.85.  Such a trade off may be valuable for problems with highly varying or discontinuous coefficients.  Also, the FSAI smoother can be easily tuned to control how aggressively its operator complexity grows.  

\section{Performance model}
\label{sec:matvec}
\begin{figure}[htb!]
\centering
\includegraphics[scale = 1.0]{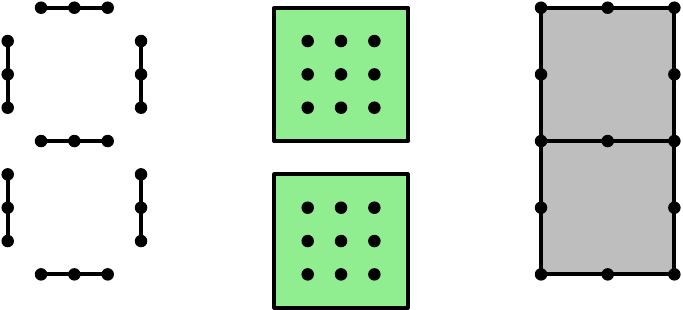}
\caption{HDG, DG, and continuous Galerkin DOFs.}
\label{hdg_fig}
\end{figure}

The storage and assembly of global matrices in finite element methods can be exceedingly prohibitive, especially at higher orders.  By leveraging matrix free algorithms, one can save on memory, and, convert a memory bound problem (sparse matrix vector multiplication) into a compute bound problem.  The authors in \cite{roca2011gpu} found that to improve sparse matrix vector multiplication for HDG methods, specialized storage formats were needed.  Since the HDG method can be reduced to a problem on the trace space, this allows for a assembly (matrix free or otherwise) of the discretization operator in a facet-by-facet manner, instead of a element-by-element manner.  The importance of this is that the element-by-element approach requires a synchronization (barrier, atomic, etc.) in order to avoid race conditions.  A graph coloring algorithm is typically used in this situation, but only allows for a group of colors to be utilized at any given time.  The facet-by-facet approach does not depend on the vertex degree, but only on the element type (how many facets on a given element).  Fig.~\ref{hdg_fig} displays some example DOFs and connectivity for HDG, DG, and continuous Galerkin methods.

The HDG method is able to exploit facet-by-facet connectivity; given an interior facet, only the two elements that share said facet will contribute to its DOFs.  An element-to-facet mapping will allow one to gather and scatter the DOFs on a facet.  Pseudocode for our matrix vector multiply routine is given in Algorithm~(\ref{algo_mv}).
\begin{center}
\begin{algorithm}
\begin{algorithmic}[1]
\Procedure{matrix vector multiply}{}
\State Algorithm for HDG matrix vector multiply (assembly free)
\State Load $\mathbb{K}_K$ and $x$ 

\State eN $\gets$ Map$0$(tid)   (Given thread id, find the edge it corresponds to)
\State (E1,E2) $\gets$ Map$1$(eN)  (Find the elements that share edge).
\State (idx$1$) $\gets$ Map$2$(eN,E1) (load local index of DOF on element 1)
\State (idx$2$) $\gets$ Map$3$(eN,E2) (load local index of DOF on element 2)

\State $y_1 = 0$

  \For{$j = 1$ to size($\mathbb{K}_\text{E1,eN}$)} 
			\State $y_1 \gets y_1 + \mathbb{K}_\text{E1,eN}($ idx$1$(j), : $) \cdot x$
      \EndFor

\State $y_2 =0$
  \For{$j = 1$ to size($\mathbb{K}_\text{E2,eN}$)}
			\State $y_2 \gets y_2 + \mathbb{K}_\text{E2,eN}($ idx$2$(j), : $)\cdot x$
      \EndFor

\State $y$(tid) $= y_1 + y_2$
\EndProcedure
\end{algorithmic}
\caption{Algorithm for HDG matrix vector multiply (assembly free).}
\label{algo_mv}
\end{algorithm}
\end{center}

%
%
%
%
%
%
%

\subsubsection{Knights Landing (KNL) many-core coprocessor}
In this work we use the second generation Xeon Phi designed by {Intel\textsuperscript{\textregistered}}.  It is a many-core processor, and the 7000 series has anywhere from 64 to 72 cores, with 4 threads per core.  The KNL is an example of a high throughput low memory device.  The design of the KNL is similar to other many-core devices and accelerators: a large number of cores with lower clock speeds are packed into the unit, enabling a large vector width, as well as having access to a user-manageable fast memory hierarchy.  One interesting feature of the KNL is that it can be programmed using traditional parallel paradigms like OpenMP, MPI, and pthreads.  In addition, it operates as a native processor.  Further details regarding the KNL can be obtained in \cite{jeffers2016intel}.  For specifications of the KNL used in this work, see Table~\ref{testbed}.
\subsubsection{Local matrix generation}
In order to make use of the assembly free matrix vector multiply, one needs to generate the associated local matrices (in the case of HDG, see (i), (iv), and (v) in Section~\ref{sec:disc0}).  To generate the required local matrices, for simplicity, we use a one core (thread) per element strategy. However, as one increases the polynomial order (beyond $p=5$), this strategy loses performance.  Further improvements may potentially be obtained by utilizing nested parallelism, or linear algebra libraries dedicated for small or medium sized matrices (\cite{heineckelibxsmm}, \cite{king2014exploiting}).

For all computational experiments we set the KNL in quadrant mode.  Figs.~\ref{fig_as1} and~\ref{fig_as2} display the wall-clock and speed up as the polynomial order and number of threads is varied.  In Fig.~\ref{fig_as1}, one can see that as we increase the number of threads beyond 32 or 64, diminishing returns are very noticeable.  To generate the results of Fig.~\ref{fig_as2}, we fix the number of threads to 64, and we vary the polynomial order.  The comparison with a serial implementation clearly shows that the additional parallelism the KNL offers is beneficial; speed ups ranging from 2X to 32X are attained.

\begin{figure}[htb!]
\centering
\hspace*{-10ex}
\begin{minipage}{0.40\textwidth}
\begin{tikzpicture}[scale = 0.9]
\begin{loglogaxis}[
grid = major,
scale only axis,
width=1.1\columnwidth,
/pgfplots/log ticks with fixed point ,
xtick={1, 2,4,8,16,32,64,128,256},
ytick={0.3,0.4,0.5,1,1.5,1.9,3,4,5,10,20},
xlabel={Number of threads},
ylabel={Time (seconds)}
]

\addplot table[x index=0,y index=1,col sep=space] {pre.dat};
\addlegendentry{$p=0$}

\addplot table[x index=0,y index=2,col sep=space] {pre.dat};
\addlegendentry{$p=1$}

\addplot table[x index=0,y index=3,col sep=space] {pre.dat};
\addlegendentry{$p=1$}

\addplot table[x index=0,y index=4,col sep=space] {pre.dat};
\addlegendentry{$p=3$}

\addplot table[x index=0,y index=5,col sep=space] {pre.dat};
\addlegendentry{$p=4$}
\end{loglogaxis}

\end{tikzpicture}
\centering
\caption{Local matrix assembly wall-clock.}
\label{fig_as1}
\end{minipage}
\hspace*{12ex}
\begin{minipage}{0.40\textwidth}
\begin{tikzpicture}[scale = 0.8]
\begin{axis}[
xlabel={Polynomial Order},
ylabel={Speedup},
,xtick=data
,ytick={6,10,17,22,28,31,35}
]
\addplot table[x index = 0, y index = 1]{sp.dat};
\end{axis}

\end{tikzpicture}
\centering
\caption{Local matrix assembly speedup.}
\label{fig_as2}
\end{minipage}
\end{figure}

\subsection{Roofline analysis}
The KNL (7210 processor number) used in this analysis runs at 2.1 GHz (double precision), a double precision processing power of of 2,199 GFLOP/s, and the STREAM memory bandwidth benchmark (triad, see \cite{stream_bench}) reports a bandwidth of 300 GB/s.  These two metrics provide the performance boundaries for arithmetic throughput and memory bandwidth limits, respectively.  For our numerical experiments, the clustering mode is set as quadrant, and cache mode is set to flat.  
\begin{table}[htb!]
\centering
 \begin{tabular}{|c|c|c|}
 \hline 
 Processor number & 7210  \\ 
 \hline 
 \# of cores & 64  \\ 
 \hline 
 Processor base frequency & 1.30 GHz \\ 
 \hline 
 Cache & 32 MB L2  \\ 
 \hline 
  RAM & 384 GB DDR4   \\ 
 \hline 
 MCDRAM & 16 GB  \\ 
 \hline 
 Instruction set & 64-bit \\ 
 \hline 
 Instruction set extension & {Intel\textsuperscript{\textregistered}} AVX-512  \\ 
 \hline 
 Operating system & CentOS 7.2  \\ 
 \hline  
 Compiler & {Intel\textsuperscript{\textregistered}} Parallel Studio XE 16 \\
  \hline 
 \end{tabular} 
\caption{KNL testbed specifications.}
\label{testbed}
\end{table}

\subsubsection{Sparse Matrix Vector Multiply (SPMV)}
In general, matrix vector multiplication is severely bandwidth bound.  For a dense matrix of dimension $N$, we can expect $(2N^2-N)$ FLOPs and $8(2N^2+N)$ MEMOPs (bytes).  Thus, for large $N$, we can expect an arithmetic intensity of only $\lim_{N\to \infty} (2N^2-N)/(8(2N^2+N))=0.25$.  For finite element problems the underlying discretization matrix is typically sparse, which pushes SPMV further into the bandwidth bound region on the roof line chart.  Matrix free operation can remedy this situation by improving the arithmetic intensity.  That is, matrix free application can shift our bandwidth bound SPMV to a compute bound (or less bandwidth bound) problem.

\begin{figure}[ht!]
\centering
    \captionsetup{justification=centering}
    \subfloat[Roofline analysis for matrix-free SPMV.]{
        \includegraphics[trim = 30mm 85mm 45mm 90mm, clip, scale = 0.45]{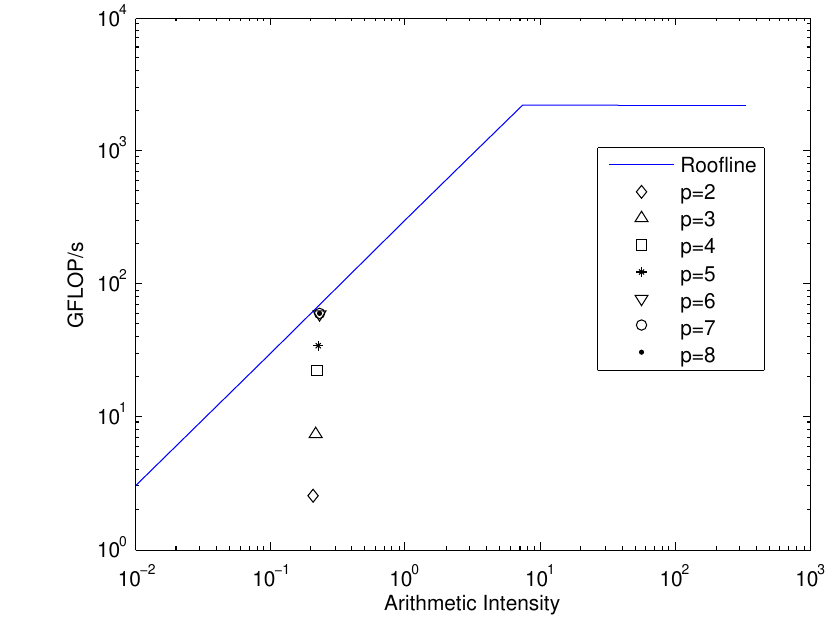}
        \label{roof1}
    }
    \subfloat[Roofline analysis for CSR SPMV.]{
        \includegraphics[trim = 30mm 85mm 45mm 90mm, clip, scale = 0.45]{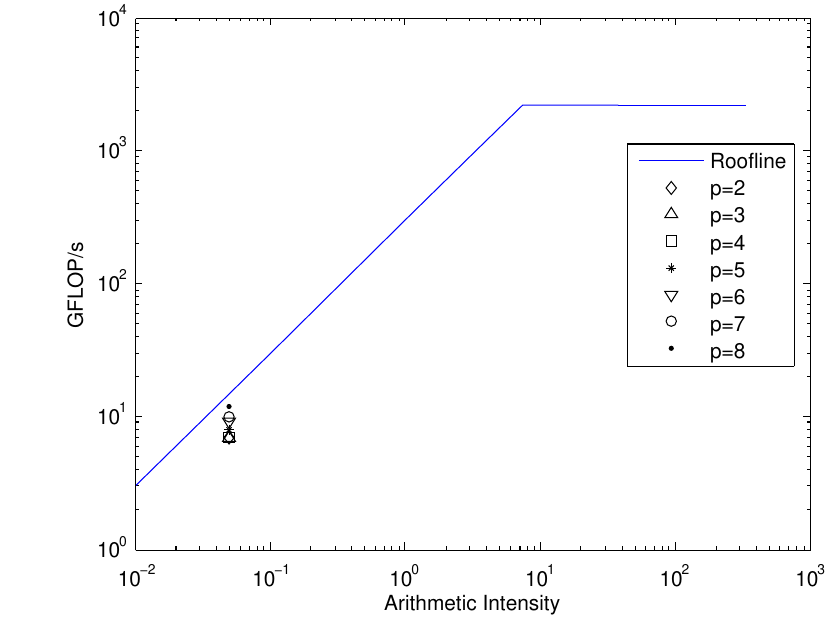}
        \label{roof2}   
    }
\captionsetup{justification=justified}
\caption{Roofline analysis for matrix-free SPMV and CSR SPMV.}
\end{figure}

In Fig.~\ref{roof1}, we include a roofline analysis of our matrix vector multiply routine; in the context of its performance on the Laplacian operator.  The roofline analysis (\cite{Williams_roof}) allows us to identify bottlenecks, confirm algorithm limitations, as well as gives us insight on what we should focus on in terms of optimization.  Two different techniques are studied.  Sparse matrix storage (CSR format, see Fig. \ref{roof2}) has the lowest performance arithmetic intensity, with a theoretical limit of 0.25.  The multithreaded Intel MKL library is used for this approach (\cite{intelMKL}).  For high orders the HDG method has hundreds of nonzeros per row.  Slightly better but similar performance was obtained in~\cite{roca2011gpu} by using a specialized block sparse matrix storage format on GPUs.  The matrix free technique shifts the arithmetic intensity favorably (see \ref{roof1}); this behavior is typical in spectral methods and spectral elements, due to near constant FLOP and memory requirements per degree of freedom (\cite{vos2010h}, \cite{cantwell2011h}, \cite{kronbichler2012generic}).

\begin{figure}[htb!]
\centering
    \captionsetup{justification=centering}
    \subfloat[Bandwidth ($p$-refinement).]{
	\includegraphics[trim = 40mm 85mm 40mm 90mm, clip, scale = 0.45]{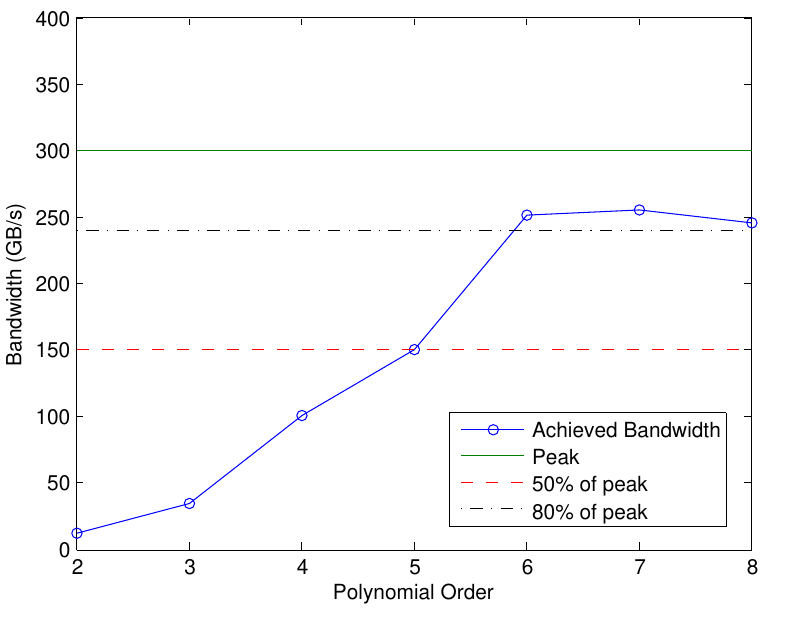}
	\label{fig:bandz1}
	}
    \subfloat[Bandwidth ($h$-refinement).]{
        \includegraphics[trim = 30mm 85mm 42mm 90mm, clip, scale = 0.45]{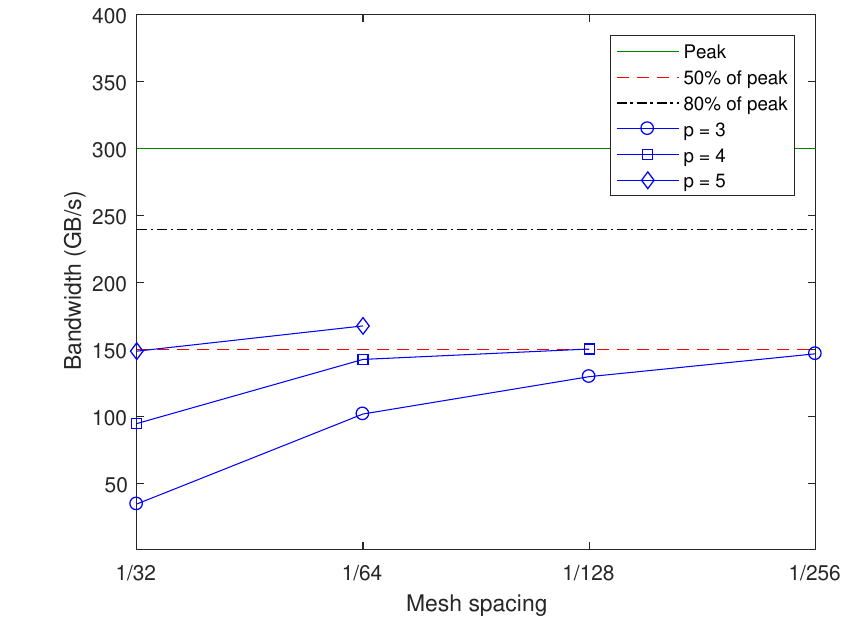}
        \label{fig:bandz2}  
    }
\captionsetup{justification=justified}
\caption{Achieved bandwidth vs polynomial order for matrix-free SPMV.}
\label{bandz1}
\end{figure}

Fig. \ref{bandz1} displays the bandwidth measurements.  From Fig.~\ref{fig:bandz1}, for lower order polynomials ($p<5$) on a fixed coarse mesh, we see that a maximum of 50\% of the peak bandwidth is achieved.  The overhead costs of forking threads is dominating computations for $p<5$ on coarser meshes.  As $p$ increases beyond five, the cost of forking threads is no longer dominant.  Our algorithm reaches 80\% of the peak bandwidth reported by the STREAM benchmark.  

Fig.~\ref{fig:bandz2} considers fixing the polynomial degree, and increasing the number of elements in the mesh.  Here it is evident that $h$-refinement improves the bandwidth performance for lower polynomial degrees.  To observe further improvements, a new SPMV routine could be designed specifically for lower orders.  For instance, in~\cite{KnepleyRT16}, the authors group together the contributions from multiple elements into a thread block in their SPMV-GPU kernel.  Their technique is shown to be very high performing.


\subsubsection{Projection From Volume To Surface}
The generation of the local solvers from equation~(\ref{shur0}) is completely data parallel and needs no synchronizations or thread communication.  It requires a dense linear solve, two matrix vector multiplications and a single SAXPY.  Roughly the complexity model for this projection operator in terms of FLOPs can be estimated by 
\begin{align*}
\text{FLOPs} &= 2N^2-N+ (2/3)N^3 + M(2N-1) + 2M), \\
\text{MEMOPs} &= 8(N^2+NM+M+N).
\end{align*}
where $N=4(p+1)^2$ and $M = 4 (k+1)$.  Table~(\ref{tabl1}) collects the results for polynomial order $0,1,2,3,$ and $4$.  As we increase the polynomial order, performance increases due to the compute bound nature of the local HDG solvers; as the dense linear solve dominates with $\mathcal{O}(n^3)$ FLOPs for $\mathcal{O}(n^2)$ data.  We again note that the focus of this work is not on the generation of this local matrices, but we can still obtain reasonable performance utilizing a straightforward implementation.

 After solving the trace space system given by equation~(\ref{shur1}), if the volume solutions $u$ and $\bm q$ are desired, one can invoke equation~(\ref{shur0}) to reconstruct them.  The cost of this procedure depends on if one discards the local matrices.  In this case, one has to recompute and the cost is the same as the projection from the surface to volume.  If one instead keeps the local matrices, all that is needed is two matrix vector multiplies and a single SAXPY.  Ultimately this comes down to a preference of convenience over memory concerns.

\def\arraystretch{1.5}
\begin{table}[htb!]
\centering
\begin{tabular}{|c|c|c|}
\hline 
Polynomial Order & AI & GFLOP/s \\ 
\hline 
0 & 0.84 & 0.0032 \\ 
\hline 
1 & 1.41 & 1.15 \\ 
\hline 
2 & 2.60 & 5.129 \\ 
\hline 
3 & 4.33 & 21.18 \\ 
\hline 
4 & 6.56 & 53.66 \\ 
\hline 
\end{tabular} 
\caption{Arithmetic intensity and GFLOP/s for the local solvers.}
\label{tabl1}
\end{table}

\subsubsection{Cost analysis matrix-free vs matrix-stored}
\begin{figure}[ht!]
\centering
    \captionsetup{justification=centering}
    \subfloat[Matrix-based solve and assembly.]{
        \includegraphics[trim = 30mm 85mm 45mm 90mm, clip, scale = 0.45]{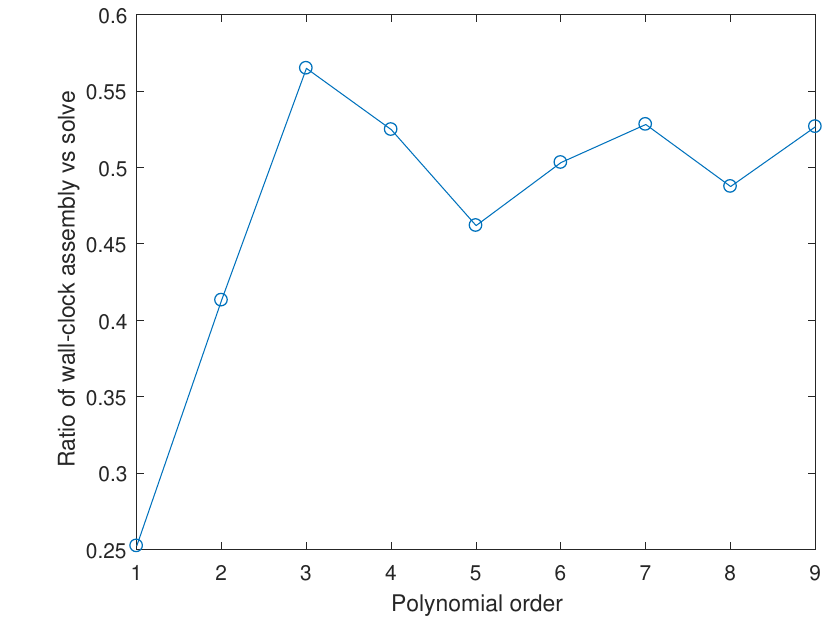}
        \label{fig:cost1}
    }
    \subfloat[Matrix-free vs matrix-based solve.]{
        \includegraphics[trim = 30mm 85mm 45mm 90mm, clip, scale = 0.45]{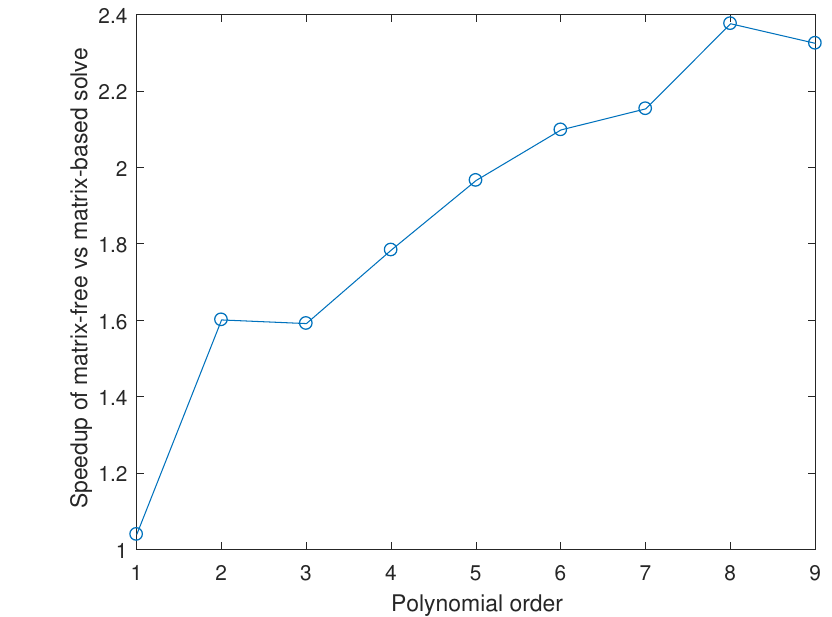}
        \label{fig:cost2}   
    }
\captionsetup{justification=justified}
\caption{Figure (a): Ratio of wall-clock timings for the matrix-based solve and assembly.  Figure(b): Speedup for matrix-free vs matrix-based solve.}
\end{figure}

We consider the cost of our matrix free approach compared to a matrix based approach.  The previous sections discuss the memory considerations, but there are other measures one can use.  Most importantly, operation counts and time-to-solution.  The operation count (FLOPs) for sparse matrix vector multiplication is $\mathcal{O}(nnz)$, where $nnz$ is the number of nonzero entries in the sparse matrix in question.  For the statically condensed HDG method on a structured Cartesian mesh, $nnz$ for the stiffness matrix is $C_1=\mathcal{O}(N d (p+1)^{d-1} \times (4d-1)(p+1)^{d-1})$ \cite{samii2016hybridized}, \cite{huerta2013efficiency} ($d$ is the dimension, $N$ is the total number of elements in the mesh).  The operation count for our proposed matrix free variant (Algorithm~\ref{algo_mv}) is $C_2=\mathcal{O}(2|\Gamma_h|( 2\cdot 4  (p+1) - 1 ))$.  In 2D, as $N\to \infty$, the ratio tends to $C_2/C_1\to(16p+14)/( 7(p+1)^2)$.  Which indicates that for higher orders, the matrix free variant is cheaper in terms of operations than the sparse matrix format.
\\
\indent  A second important measure is time-to-solution.  Both the matrix free and matrix based algorithm require the local solvers.  In the matrix based approach, a non negligible amount of time is spent on the matrix assembly, especially for higher orders.  We compare the wall-clock matrix based solve to the wall-clock of the (parallel) assembly.  See Fig.~\ref{fig:cost1}.  For $p=2$, the assembly time is already more than 40\% of the total solve time.  Increasing the polynomial order to $p>2$, the assembly time ranges from around 45\% to 56\% of the total solve time.  Comparing the matrix free solve timings, to the matrix based solve and assembly time, Fig.~\ref{fig:cost2} shows for $p>2$ we obtain speedups of 1.5 to 2.3 in favor of the matrix free method.

We acknowledge that the benefits of matrix free methods will vanish if efficient preconditioners are not available.  In such cases, in order to have more options for robust preconditioners, the matrix will have to be stored.

\section{Conclusions}
We presented a highly efficient multigrid technique for solving a second order elliptic PDE with the HDG discretization.  Both the multigrid technique and the HDG discretization are amenable to high throughput low memory environments like that provided by the KNL architecture; this was demonstrated by using thorough profiling and utilizing a roofline analysis (\cite{Williams_roof}).  The HDG method has much of its computation localized due to the discontinuous nature of the solution technique.  Moreover, the HDG method brings static condensation to DG methods, which significantly reduces the number of nonzeros in the discretization operator.  This in turn means that less work is required for a linear solver to obtain solutions.  Since the HDG method converges with optimal orders for all of the variables it approximates, a local element by element postprocessing is available.  

Any iterative solver will require sparse matrix vector multiplication (matrix free or otherwise); including multigrid.  Two different approaches were examined for sparse matrix vector multiplication: matrix based and matrix free.  For high order HDG, the matrix free technique better utilizes the resources available on the KNL, because of increased arithmetic intensity.  In the high order regime, the matrix based technique requires a sparse matrix storage which increases time to solution.  This is mainly because of additional assembly time, having a low arithmetic intensity, and erratic memory access patterns for sparse matrix vector multiplication.

Our algorithm is able to attain 80\% of peak bandwidth performance for higher order polynomials.  This is possible due to the data locality inherent in the HDG method. Very good performance is realized for high order schemes, due to good arithmetic intensity, which declines as the order is reduced.  The performance is lower for polynomial orders $p<6$.  With an approach similar to the work done in \cite{KnepleyRT16}, where multiple cells are processed by a thread block, one would be able to increase performance for lower order polynomials.  We observed speedups when compared to a multicore CPU for the HDG methods components, namely, volume to surface mapping, surface to volume mapping, matrix free matrix vector multiplication, and local matrix generation.  The ratio of peak flop rates on the two target architectures was roughly 100X, and peak bandwidths was 5X, so this figure fits with our model of the computation.  This is possible due to the data locality inherent in the HDG method. Very good performance is realized for high order schemes, due to good arithmetic intensity, which declines as the order is reduced.
    
An attractive feature of the KNL is that it can be programmed with traditional parallel paradigms like OpenMP, MPI, and pthreads.  Thus, one can harness the massive fine-grain parallelism that the KNL offers by utilizing these traditional parallel paradigms with significantly limited intrusion.

\section*{Acknowledgments}
Fabien acknowledges the support from the Ken Kennedy Institute and the Ken Kennedy--Cray graduate Fellowship.  Fabien also acknowledges the support from the Richard Tapia Center for Excellence \& Equity, and the Rice Graduate Education for Minorities program.  This work used the Extreme Science and Engineering Discovery Environment (XSEDE), which is supported by National Science Foundation grant number ACI-1053575.  The KNL used in this work was provided to the RELACS research group by Richard T. Mills through the {Intel\textsuperscript{\textregistered}} Parallel Computing Center at Rice University.  

\bibliographystyle{siamplain}
\bibliography{references}

\begin{thebibliography}{10}

\bibitem{Adams2003593}
{\sc M.~Adams, M.~Brezina, J.~Hu, and R.~Tuminaro}, {\em Parallel multigrid
  smoothing: polynomial versus {G}auss--{S}eidel}, Journal of Computational
  Physics, 188 (2003), pp.~593 -- 610.

\bibitem{Anzt2012}
{\sc H.~Ali, Y.~Shi, D.~Khazanchi, M.~Lees, G.~D. van Albada, J.~Dongarra,
  P.~M. Sloot, J.~Dongarra, H.~Anzt, S.~Tomov, M.~Gates, J.~Dongarra, and
  V.~Heuveline}, {\em Proceedings of the {I}nternational {C}onference on
  {C}omputational {S}cience, {ICCS} 2012 {B}lock-asynchronous multigrid
  smoothers for {GPU}-accelerated systems}, Procedia Computer Science, 9
  (2012), pp.~7 -- 16.

\bibitem{Ant}
{\sc P.~F. Antonietti, M.~Sarti, and M.~Verani}, {\em Multigrid algorithms for
  hp-discontinuous {G}alerkin discretizations of elliptic problems}, SIAM
  Journal on Numerical Analysis,  (2015).

\bibitem{Antonietti2016}
{\sc P.~F. Antonietti, M.~Sarti, and M.~Verani}, {\em Multigrid algorithms for
  high order discontinuous {G}alerkin methods}, Domain Decomposition Methods in
  Science and Engineering XXII,  (2016), pp.~3--13.

\bibitem{Anzt2013}
{\sc H.~Anzt, S.~Tomov, J.~Dongarra, and V.~Heuveline}, {\em Weighted
  block-asynchronous iteration on {GPU}-accelerated systems},  (2013),
  pp.~145--154.

\bibitem{ArnoldBCM02}
{\sc D.~N. Arnold, B.~C. F.~Brezzi, and L.~D. Marini}, {\em Unified analysis of
  discontinuous {G}alerkin methods for elliptic problems}, {SIAM} J. Numerical
  Analysis, 39 (2002), pp.~1749--1779.

\bibitem{Bastian}
{\sc P.~Bastian}, {\em Load balancing for adaptive multigrid methods}, SIAM
  Journal on Scientific Computing, 19 (1998), pp.~1303--1321.

\bibitem{bastian1998additive}
{\sc P.~Bastian, G.~Wittum, and W.~Hackbusch}, {\em Additive and multiplicative
  multi-grid; a comparison}, Computing, 60 (1998), pp.~345--364.

\bibitem{Benson}
{\sc M.~W. Benson}, {\em Frequency domain behavior of a set of parallel
  multigrid smoothing operators}, International Journal of Computer
  Mathematics, 36 (1990), pp.~77--88.

\bibitem{BerrutT04}
{\sc J.~P. Berrut and L.~N. Trefethen}, {\em Barycentric {L}agrange
  interpolation}, {SIAM} Review, 46 (2004), pp.~501--517.

\bibitem{Brandt1}
{\sc A.~Brandt}, {\em Multi-level adaptive solutions to boundary-value
  problems}, Mathematics of Computation, 31 (1977), pp.~333--390.

\bibitem{Brandt2}
{\sc A.~Brandt and O.~Livne}, {\em Multigrid Techniques: 1984 Guide with
  Applications to Fluid Dynamics, Revised Edition}, Classics in Applied
  Mathematics, Society for Industrial and Applied Mathematics, 2011.

\bibitem{brenner2009multigrid}
{\sc S.~C. Brenner, J.~Cui, and L.-Y. Sung}, {\em Multigrid methods for the
  symmetric interior penalty method on graded meshes}, Numerical Linear Algebra
  with Applications, 16 (2009), pp.~481--501.

\bibitem{brenner2005convergence}
{\sc S.~C. Brenner and J.~Zhao}, {\em Convergence of multigrid algorithms for
  interior penalty methods}, Applied Numerical Analysis \& Computational
  Mathematics, 2 (2005), pp.~3--18.

\bibitem{Broker}
{\sc O.~Br\"{o}ker and M.~J. Grote}, {\em Sparse approximate inverse smoothers
  for geometric and algebraic multigrid}, Appl. Numer. Math., 41 (2002),
  pp.~61--80.

\bibitem{cantwell2011h}
{\sc C.~Cantwell, S.~Sherwin, R.~Kirby, and P.~Kelly}, {\em From h to p
  efficiently: Strategy selection for operator evaluation on hexahedral and
  tetrahedral elements}, Computers \& Fluids, 43 (2011), pp.~23--28.

\bibitem{CockburnDGRS09}
{\sc B.~Cockburn, J.~G. B.~Dong, M.~Restelli, and R.~Sacco}, {\em A
  hybridizable discontinuous {G}alerkin method for steady-state
  convection-diffusion-reaction problems}, {SIAM} J. Scientific Computing, 31
  (2009), pp.~3827--3846.

\bibitem{CockburnDG08}
{\sc B.~Cockburn, B.~Dong, and J.~Guzm{\'{a}}n}, {\em A superconvergent
  {LDG}-hybridizable {G}alerkin method for second-order elliptic problems},
  Math. Comput., 77 (2008), pp.~1887--1916.

\bibitem{cockburn2013multigrid}
{\sc B.~Cockburn, O.~Dubois, J.~Gopalakrishnan, and S.~Tan}, {\em Multigrid for
  an {HDG} method}, IMA Journal of Numerical Analysis, 34 (2013),
  pp.~1386--1425.

\bibitem{CockburnGL09}
{\sc B.~Cockburn, J.~Gopalakrishnan, and R.~D. Lazarov}, {\em Unified
  hybridization of discontinuous {G}alerkin, mixed, and continuous {G}alerkin
  methods for second order elliptic problems}, {SIAM} J. Numerical Analysis, 47
  (2009), pp.~1319--1365.

\bibitem{Fedorenko}
{\sc R.~P. Fedorenko}, {\em The speed of convergence of one iterative process},
  U.S.S.R. Comput. Math. Math. Phys., 4 (1964), pp.~559--564.

\bibitem{Fidkowski}
{\sc K.~J. Fidkowski, T.~A. Oliver, J.~Lu, and D.~L. Darmofal}, {\em
  p-multigrid solution of high-order discontinuous {G}alerkin discretizations
  of the compressible navier-stokes equations}, J. Comput. Phys., 207 (2005),
  pp.~92--113.

\bibitem{gopalakrishnan2003multilevel}
{\sc J.~Gopalakrishnan and G.~Kanschat}, {\em A multilevel discontinuous
  {G}alerkin method}, Numerische Mathematik, 95 (2003), pp.~527--550.

\bibitem{gopalakrishnan2009convergent}
{\sc J.~Gopalakrishnan and S.~Tan}, {\em A convergent multigrid cycle for the
  hybridized mixed method}, Numerical Linear Algebra with Applications, 16
  (2009), pp.~689--714.

\bibitem{hackbusch}
{\sc W.~Hackbusch}, {\em Multi-Grid Methods and Applications}, Springer Series
  in Computational Mathematics, Springer Berlin Heidelberg, 2013.

\bibitem{heineckelibxsmm}
{\sc A.~Heinecke, H.~Pabst, and G.~Henry}, {\em {LIBXSMM}: A high performance
  library for small matrix multiplications}.

\bibitem{hemker2003two}
{\sc P.~Hemker, W.~Hoffmann, and M.~Van~Raalte}, {\em Two-level fourier
  analysis of a multigrid approach for discontinuous {G}alerkin
  discretization}, SIAM Journal on Scientific Computing, 25 (2003),
  pp.~1018--1041.

\bibitem{hemker2004fourierlinear}
{\sc P.~Hemker, W.~Hoffmann, and M.~Van~Raalte}, {\em Fourier two-level
  analysis for discontinuous {G}alerkin discretization with linear elements},
  Numerical linear algebra with applications, 11 (2004), pp.~473--491.

\bibitem{hemker2004fourier}
{\sc P.~Hemker and M.~Van~Raalte}, {\em Fourier two-level analysis for higher
  dimensional discontinuous {G}alerkin discretisation}, Computing and
  Visualization in Science, 7 (2004), pp.~159--172.

\bibitem{huerta2013efficiency}
{\sc A.~Huerta, A.~Angeloski, X.~Roca, and J.~Peraire}, {\em Efficiency of
  high-order elements for continuous and discontinuous {G}alerkin methods},
  International Journal for Numerical Methods in Engineering, 96 (2013),
  pp.~529--560.

\bibitem{intelMKL}
{\em Intel {M}ath {K}ernel {L}ibrary. reference manual}, 2009.

\bibitem{jeffers2016intel}
{\sc J.~Jeffers, J.~Reinders, and A.~Sodani}, {\em Intel Xeon Phi Processor
  High Performance Programming: Knights Landing Edition}, Morgan Kaufmann,
  2016.

\bibitem{johannsen2005multigrid}
{\sc K.~Johannsen}, {\em Multigrid methods for {NIPG}}, ICES Report,  (2005),
  pp.~05--31.

\bibitem{johannsen2005symmetric}
{\sc K.~Johannsen}, {\em A symmetric smoother for the nonsymmetric interior
  penalty discontinuous {G}alerkin discretization}, ICES Report, 5 (2005),
  p.~23.

\bibitem{karniadakis1999spectral}
{\sc G.~Karniadakis and S.~J. Sherwin}, {\em Spectral/hp Element Methods for
  CFD}, Numerical mathematics and scientific computation, Oxford University
  Press, 1999.

\bibitem{king2014exploiting}
{\sc J.~King, S.~Yakovlev, Z.~Fu, R.~M. Kirby, and S.~J. Sherwin}, {\em
  Exploiting batch processing on streaming architectures to solve {2D} elliptic
  finite element problems: A hybridized discontinuous {G}alerkin ({HDG}) case
  study}, Journal of Scientific Computing, 60 (2014), pp.~457--482.

\bibitem{Kirby2012}
{\sc R.~M. Kirby, S.~J. Sherwin, and B.~Cockburn}, {\em To {CG} or to {HDG}: A
  comparative study}, Journal of Scientific Computing, 51 (2012), pp.~183--212.

\bibitem{KnepleyRT16}
{\sc M.~G. Knepley, K.~Rupp, and A.~R. Terrel}, {\em Finite element integration
  with quadrature on the {GPU}}, CoRR, abs/1607.04245 (2016).

\bibitem{kronbichler2012generic}
{\sc M.~Kronbichler and K.~Kormann}, {\em A generic interface for parallel
  cell-based finite element operator application}, Computers \& Fluids, 63
  (2012), pp.~135--147.

\bibitem{Loisel2008}
{\sc S.~Loisel, R.~Nabben, and D.~B. Szyld}, {\em On hybrid multigrid-{S}chwarz
  algorithms}, Journal of Scientific Computing, 36 (2008), pp.~165--175.

\bibitem{Lottes2005}
{\sc J.~W. Lottes and P.~F. Fischer}, {\em Hybrid multigrid/{S}chwarz
  algorithms for the spectral element method}, Journal of Scientific Computing,
  24 (2005), pp.~45--78.

\bibitem{stream_bench}
{\sc J.~D. McCalpin}, {\em {STREAM}: Sustainable memory bandwidth in high
  performance computers}.
\newblock https://www.cs.virginia.edu/stream/.
\newblock Accessed: 2016-08-30.

\bibitem{Mitchell}
{\sc W.~F. Mitchell}, {\em The hp-multigrid method applied to hp-adaptive
  refinement of triangular grids}, Numerical Linear Algebra with Applications,
  17 (2010), pp.~211--228.

\bibitem{napov2010does}
{\sc A.~Napov and Y.~Notay}, {\em When does two-grid optimality carry over to
  the v-cycle?}, Numerical linear algebra with applications, 17 (2010),
  pp.~273--290.

\bibitem{roca2011gpu}
{\sc X.~Roca, N.~C. Nguyen, and J.~Peraire}, {\em {GPU}-accelerated sparse
  matrix-vector product for a hybridizable discontinuous {G}alerkin method}, in
  Aerospace Sciences Meetings. American Institute of Aeronautics and
  Astronautics, 2011, pp.~2011--687.

\bibitem{saad}
{\sc Y.~Saad}, {\em Iterative Methods for Sparse Linear Systems: Second
  Edition}, Society for Industrial and Applied Mathematics, 2003.

\bibitem{samii2016hybridized}
{\sc A.~Samii, N.~Panda, C.~Michoski, and C.~Dawson}, {\em A hybridized
  discontinuous {G}alerkin method for the nonlinear {K}orteweg--de {V}ries
  equation}, Journal of Scientific Computing, 68 (2016), pp.~191--212.

\bibitem{schenk2004solving}
{\sc O.~Schenk and K.~G{\"a}rtner}, {\em Solving unsymmetric sparse systems of
  linear equations with {PARDISO}}, Future Generation Computer Systems, 20
  (2004), pp.~475--487.

\bibitem{Stiller15}
{\sc J.~Stiller}, {\em Nonuniformly weighted {S}chwarz smoothers for spectral
  element multigrid}, CoRR, abs/1512.02390 (2015).

\bibitem{Stiller16}
{\sc J.~Stiller}, {\em Robust multigrid for high-order discontinuous {G}alerkin
  methods: {A} fast {P}oisson solver suitable for high-aspect ratio cartesian
  grids}, CoRR, abs/1603.02524 (2016).

\bibitem{tan09}
{\sc S.~Tan}, {\em Iterative solvers for hybridized finite element methods},
  {PhD} Thesis,  (2009).

\bibitem{trefethen2013approximation}
{\sc L.~N. Trefethen}, {\em Approximation Theory and Approximation Practice},
  Siam, 2013.

\bibitem{trottenberg}
{\sc U.~Trottenberg, C.~W. Oosterlee, and A.~Sch{\"u}ller}, {\em Multigrid},
  Academic Press, 2001.

\bibitem{vos2010h}
{\sc P.~E. Vos, S.~J. Sherwin, and R.~M. Kirby}, {\em From h to p efficiently:
  Implementing finite and spectral/hp element methods to achieve optimal
  performance for low-and high-order discretisations}, Journal of Computational
  Physics, 229 (2010), pp.~5161--5181.

\bibitem{wang2014explicit}
{\sc H.~Wang, D.~Huybrechs, and S.~Vandewalle}, {\em Explicit barycentric
  weights for polynomial interpolation in the roots or extrema of classical
  orthogonal polynomials}, Mathematics of Computation, 83 (2014),
  pp.~2893--2914.

\bibitem{wesseling}
{\sc P.~Wesseling}, {\em An introduction to multigrid methods}, Pure and
  applied mathematics, John Wiley \& Sons Australia, Limited, 1992.

\bibitem{Williams_roof}
{\sc S.~Williams, A.~Waterman, and D.~Patterson}, {\em Roofline: An insightful
  visual performance model for multicore architectures}, Commun. ACM, 52
  (2009), pp.~65--76.

\bibitem{yakovlev2016cg}
{\sc S.~Yakovlev, D.~Moxey, R.~M. Kirby, and S.~J. Sherwin}, {\em To {CG} or to
  {HDG}: A comparative study in {3D}}, Journal of Scientific Computing, 67
  (2016), pp.~192--220.

\end{thebibliography}
\end{document}